\documentclass[11pt]{amsart}
\usepackage{amsbsy,amssymb,amscd,amsfonts,latexsym,amstext,delarray,
amsmath,amsthm, graphicx}

\input xypic

\newtheorem{thm}{Theorem}[section]
\newtheorem{prop}[thm]{Proposition}
\newtheorem{cor}[thm]{Corollary}
\newtheorem{lem}[thm]{Lemma}

\newtheorem{defn}[thm]{Definition}

\numberwithin{equation}{section}

\def\bM{{\mathbb M}}

\def\bT{{\mathbb T}}

\def\C{{\mathbb C}}

\renewcommand{\H}{{\mathbb H}}

\renewcommand{\P}{{\mathbb P}}
\def\Q{{\mathbb Q}}
\def\R{{\mathbb R}}
\def\Z{{\mathbb Z}}
\def\K{{\mathbb K}}

\def\cA{{\mathcal A}}
\def\cB{{\mathcal B}}

\def\cD{{\mathcal D}}

\def\cF{{\mathcal F}}

\def\cH{{\mathcal H}}

\def\cK{{\mathcal K}}

\def\cR{{\mathcal R}}

\def\cU{{\mathcal U}}
\def\cV{{\mathcal V}}

\def\Aut{{\rm Aut}}
\def\Coker{{\rm Coker}}

\def\End{{\rm End}}

\def\GL{{\rm GL}}

\def\Hom{{\rm Hom}}

\def\Ind{{\rm Ind}}

\def\Ker{{\rm Ker}}

\def\Res{{\rm Res}}
\def\sign{{\rm sign}}
\def\SL{{\rm SL}}
\def\Spec{{\rm Spec}}

\def\Tr{{\rm Tr}}
\def\tr{{\rm tr}}

\def\cancel#1#2{\ooalign{$\hfil#1\mkern1mu/\hfil$\crcr$#1#2$}}
\def\dirac{\mathpalette\cancel\partial}
\def\Dirac{\mathpalette\cancel D}
\newcommand{\ie}{{\it i.e.\/}\ }
\newcommand{\eg}{{\it e.g.\/}\ }
\newcommand{\cf}{{\it cf.\/}\ }

\title[Solvmanifolds and NC tori]{Solvmanifolds and 
noncommutative tori with real
multiplication}

\author[Marcolli]{Matilde Marcolli}
\address{Max--Planck Institut f\"ur Mathematik  \\
Vivatsgasse 7 \\ Bonn, D-53111 Germany}
\email{marcolli\@@mpim-bonn.mpg.de}

\begin{document}
\maketitle

\begin{abstract}
We prove that the Shimizu $L$-function of a real quadratic field is obtained from
a (Lorentzian) spectral triple on a noncommutative torus with real multiplication,
as an adiabatic limit of the Dirac operator on a 3-dimensional solvmanifold. 
The Dirac operator on this 3-dimensional geometry gives, via the Connes--Landi 
isospectral deformations, a spectral triple for the noncommutative tori obtained by
deforming the fiber tori to noncommutative spaces.
The 3-dimensional solvmanifold is the homotopy quotient in the sense of Baum--Connes 
of the noncommutative space obtained as the crossed product of the noncommutative 
torus by the action of the units of the real quadratic field. 
This noncommutative space is identified
with the twisted group $C^*$-algebra of the fundamental group of the 3-manifold.
The twisting can be interpreted as the cocycle arising from a magnetic field,
as in the theory of the quantum Hall effect. We prove a twisted index theorem that
computes the range of the trace on the $K$-theory of this noncommutative space
and gives an estimate on the gaps in the spectrum of the associated Harper operator.
\end{abstract}

\tableofcontents

\section{Introduction}

In the 1970s Hirzebruch formulated a conjecture, \cf \cite{Hirz}, on the topological interpretation
of certain special values of $L$-functions of totally real fields in terms of signature defects.
The conjecture was proved in the early '80s by Atiyah--Donnelly--Singer \cite{ADS} and 
by M\"uller \cite{Muller}. Hirzebruch's conjecture played an important role in the development of
the Atiyah--Patodi--Singer index theorem \cite{APS}, which in turn is a key ingredient in the proof
\cite{ADS} of the conjecture, extending the Hirzebruch--Riemann--Roch theorem to manifolds 
with boundary and relating the signature defect to the eta invariant. 
Geometrically, the link of an isolated 
singularity of the Hilbert modular variety associated to a totally real number field is given by
a $(4k-1)$-dimensional solvmanifold. The signature of the Hilbert modular variety
is then computed by the APS theorem applied to the resulting manifold with boundary and the
signature defects are computed by the eta invariant of the solvmanifold. 
The main step in the proof of \cite{ADS} then consists
of separating out the eta function of the signature operator on the solvmanifold into a part that
recovers the Shimizu $L$-function of the totally real field and a residual part, which is shown 
not to contribute to the eta invariant.

\medskip

We concentrate here on the simplest case, that or real quadratic fields, and we consider the 
question of whether the Shimizu $L$-function can be related in a similar way to a spectral 
geometry (in the sense of Connes' notion of spectral triples \cite{Co-S3}) on a noncommutative 
torus with real multiplication by the same real quadratic field.

\medskip

The noncommutative tori associated to quadratic
irrationalities have been extensively studied by Manin in
\cite{Man} and subsequently by several authors.
They have the special property of
``real multiplication'', derived from the presence of
non-trivial self Morita equivalences. 
It is argued in \cite{Man} that the noncommutative tori with real 
multiplication should play a role for real quadratic field parallel 
to the the theory of elliptic curves with complex multiplication in 
the case of imaginary quadratic fields. This makes it an interesting
problem to study the geometric properties of this particular class
of noncommutative spaces, and their relation to the arithmetic of
real quadratic fields.

\medskip

We show in \S \ref{SconjBC} that there is a close relation between 
the 3-dimensional solvmanifold and the noncommutative torus with real multiplication. 
Namely, we prove that the first is the homotopy quotient, in 
the sense of Baum--Connes, of the noncommutative space given by the quotient of the 
latter by the action of the infinite group of units, \cf \S \ref{Snctori}.
We also show that the 3-manifold can be identified with the pullback of the universal
family of elliptic curves along a closed geodesic in the modular curve.

\smallskip

This in terpretation as the homotopy quotient of a noncommutative space
provides a geometric setting analogous to the one developed in the noncommutative 
geometry models of the quantum Hall effect \cite{Bel1}, where the presence of a magnetic field makes 
the Brillouin zone of the lattice into a noncommutative torus. Here, the 3-dimensional 
solvmanifold is similarly related to a noncommutative space whose algebra of coordinates 
is the crossed product of the algebra of the noncommutative torus by the action of the units. 
This is obtained by twisting the group ring of the fundamental group of the solvmanifold by a 
cocycle, defined in terms of a magnetic potential. The noncommutative space is the resulting
twisted group $C^*$-algebra, \cf \S \ref{Stwist}. 
As in the case of the quantum Hall effect, and in the 
noncommutative Bloch theory of electron-ion interactions, one obtains in \S \ref{SrangeTr} 
information on the spectral theory of the corresponding magnetic Laplacian by computing
the range of the trace on the $K$-theory of the twisted group $C^*$-algebra. We prove a
twisted index theorem which we use to compute the range of the trace using a spectral
flow computation and the Baum--Connes conjecture, which is 
known to hold for the fundamental group of the 3-dimensional solvmanifold.

\smallskip

This way of passing from the 3-dimensional solvmanifold to the associated noncommutative space 
is obtained in two steps. Viewing the 3-manifold as a fibration of 2-dimensional tori 
over the circle, one first replaces the fiber tori by noncommutative tori and then the
mapping torus by the dual action of the units. We show in \S \ref{Sisospec}
that the first step can be seen as a case of the 
Connes--Landi isospectral deformations \cite{CoLa}. In particular, we prove that the Dirac 
operator on the 3-manifold induces in this way a Dirac operator on the noncommutative torus 
with real multiplication. A unitary equivalence as the one considered in \cite{ADS} then 
factors this Dirac operator into a product of two operators, one of which has spectrum 
given by the norms $N(\lambda)$ of the lattice points $\lambda$ and recovers the 
Shimuzu $L$-function. In \S \ref{DiffNCtori} we show how an adiabatic limit relates the Dirac 
operator on the 3-manifold
to known differential operators on the noncommutative torus, obtained by considering the
derivations along the leaves of the Kronecker foliations associated to the Galois
conjugate elements $\theta$ and $\theta'$ in the real quadratic field $\K=\Q(\theta)$. 

\smallskip

Finally we prove in \S \ref{Slorentz} that the norms $N(\lambda)$ define
the momenta of a {\em Lorentzian} Dirac operator on the noncommutative 
torus with real multiplication. The theory 
of spectral triples in Lorentzian signature is at present still under active development and this
provides a natural example where the arithmetic structure dictates how the Lorentzian
geometry should be treated in the noncommutative context. We develop a framework for
Lorentzian spectral triples over real quadratic fields, where the Galois involution of $\K$
provides a canonical choice of a Krein involution. In particular, we prove that, in passing
from the indefinite Lorentzian geometry defined by the quadratic form given by the norm to
the associated real Hilbert space, one can resolve the infinite multiplicities in the
spectrum of the Dirac operator arising from the presence of a non-compact group of
symmetries (the units of the real quadratic field). We show that the resulting operator 
on the real Hilbert space has the same eta function as the one coming from the adiabatic 
limit of the 3-dimensional geometry as in \S \ref{DiffNCtori}. 
This eta function is the product of the Shimizu $L$-function of the real quadratic 
field by a term that only depends on the fundamental unit.

\bigskip

{\bf Acknowledgment.} I am very grateful to Sir Michael Atiyah for asking 
the question this paper is attempting to answer. I thank Alain Connes, Yuri
Manin, and Don Zagier for useful conversations. 
I also thank the Mittag Leffler Institute for hospitality and support, 
while part of this work was done.
This research is partially supported by NSF grant DMS-0651925. 

\section{3-dimensional solvmanifolds and real quadratic fields}\label{S3dimsolv}

Let $\K=\Q(\sqrt{d})$ be a real quadratic field and let 
$\iota_i :\K \hookrightarrow \R$, for $i=1,2$, be its
two real embeddings. We let $L\subset \K$ be a lattice, 
with $U_L^+$ the group of totally positive units preserving $L$,
\begin{equation}\label{ULambdaplus}
U_L^+ =\{ u\in O_\K^* \,|\, u L \subset L, \ 
\iota_i(u)\in \R^*_+ \}.
\end{equation}
We denote by $\epsilon_L$ a generator, so that 
$U_L^+=\epsilon_L^\Z$. In the case where $L=O_\K$,
the ring of integers of $\K$, then the generator 
$\epsilon=\epsilon_L$ is a fundamental unit.
We consider the embedding of $L$ in $\R^2$ given
by the mapping
\begin{equation}\label{R2lattice}
L \ni \ell \mapsto (\iota_1(\ell),\iota_2(\ell))\subset \R^2.
\end{equation}
We denote the range by $\Lambda=(\iota_1,\iota_2)(L)$. 
This is a lattice in $\R^2$.
The action of $U^+_L$ extends to an action on $\Lambda$ by
\begin{equation}\label{Uplusact}
\lambda=(\iota_1(\ell),\iota_2(\ell))\mapsto (\epsilon \iota_1(\ell),
\epsilon' \iota_2(\ell)) = (\epsilon \iota_1(\ell),
\epsilon^{-1} \iota_2(\ell)).
\end{equation}

\subsection{Semidirect products and solvmanifolds}

Let us denote by $V$ either the group $U_L^+$ or a finite index 
subgroup thereof. 
As in \cite{ADS}, we consider the crossed product
\begin{equation}\label{crossV}
S(\Lambda,V)=\Lambda \rtimes_\epsilon V ,
\end{equation}
where the action of $V=\epsilon^\Z$ on $\Lambda$ is induced by the action by 
multiplication on $L$.
As shown in \cite{ADS}, these are discrete subgroups of the solvable 
Lie group
\begin{equation}\label{solvLie}
S(\R^2,\R)=\R^2\rtimes \R,
\end{equation}
with the action of $\R$ on $\R^2$ by the one parameter subgroup 
$\Theta_t(x,y)=(e^t x,e^{-t} y)$ of $\SL_2(\R)$.
For $\epsilon >1$ and $\epsilon'=\epsilon^{-1}<1$, 
the action of $V$ on $\Lambda$ is then generated by 
\begin{equation}\label{Aepsilon}
 A_\epsilon = \left(\begin{matrix} 
\epsilon & 0 \\ 0 & \epsilon' \end{matrix}\right) \in \SL_2(\R).
\end{equation}

We also consider the 3-dimensional solvmanifold obtained as the quotient
\begin{equation}\label{Xmanifold}
X_\epsilon =S(\Lambda,V)\backslash S(\R^2,\R),
\end{equation}
with $\pi_1(X_\epsilon)=S(\Lambda,V)$. 

\subsection{The topology of the 3-manifold $X_\epsilon$}\label{S3dimtop}

It is well known \cite{ADS} that the 3-manifold $X_\epsilon$ of
\eqref{Xmanifold} is a fibration over the circle $S^1$ with fibers
that are 2-tori and with monodromy given by the matrix $A_\epsilon$ of
\eqref{Aepsilon}.

\begin{lem}\label{H1Xlem}
The manifold $X_\epsilon$ has first homology
\begin{equation}\label{H1Xe}
H_1(X_\epsilon,\Z) = \Lambda/(1-A_\epsilon)\Lambda \oplus \Z.
\end{equation}
\end{lem}

\proof The fundamental group is $\pi_1(X_\epsilon)=S(\Lambda,V)$.
Consider the surjective map
\begin{equation}\label{piS}
 \pi: S(\Lambda,V) \to \Lambda/(1-A_\epsilon)\Lambda \oplus \Z ,  \ \
\  \pi(\lambda,n)=(\lambda \mod (1-A_\epsilon)\Lambda ,n).
\end{equation}
By writing
$$ A_\epsilon^n (\lambda')=\lambda' - (1-A_\epsilon)(\lambda'+A_\epsilon
(\lambda') +\cdots+ A_\epsilon^{n-1}(\lambda')) $$
one sees that $\lambda +A_\epsilon^n(\lambda') =
\lambda + \lambda'$ modulo $(1-A_\epsilon)\Lambda$, so that 
$\pi(\lambda + A_\epsilon^n  \lambda', n+n')=\pi(\lambda,n)+
\pi(\lambda',n')$. Since commutators in $S(\Lambda,V)$ are of the form
$$ (\lambda,n)(\lambda',n')(\lambda,n)^{-1}(\lambda',n')^{-1} =
((1-A_\epsilon^{n'})\lambda -(1-A_\epsilon^n)\lambda',0), $$
we see that the homomorphism \eqref{piS} has 
$\Ker(\pi)=[\pi_1(X_\epsilon),\pi_1(X_\epsilon)]$. 
\endproof

\begin{cor}\label{cohomXe}
The compact 3-manifold $X_\epsilon$ has cohomology
\begin{equation}\label{Hevodd}
H^{even}(X_\epsilon,\Z)=\Z \oplus \Z \oplus
\Coker(1-A_\epsilon), \ \ \ 
H^{odd}(X_\epsilon,\Z)=\Z \oplus \Z. 
\end{equation}
\end{cor}

\proof
By Poincar\'e duality we have
\begin{equation}\label{Poincdual}
H^2(X_\epsilon,\Z)\cong H_1(X_\epsilon,\Z)\cong \Z \oplus
\Lambda/(1-A_\epsilon)\Lambda,
\end{equation}
and $H^1(X_\epsilon,\Z)=\Hom(H_1(X_\epsilon,\Z),\Z)=\Z \oplus
\Hom(\Lambda/(1-A_\epsilon),\Z)$, so that 
\begin{equation}\label{Poincdual2}
H_2(X_\epsilon,\Z)\cong H^1(X_\epsilon,\Z)\cong \Z \oplus
\Hom(\Lambda/(1-A_\epsilon)\Lambda,\Z).
\end{equation}
We have $\Ker(1-A_\epsilon)=0$, while $\Coker(1-A_\epsilon)$ is torsion,
so that we obtain
\begin{equation}\label{H1H2}
\begin{array}{cl}
H^0(X_\epsilon,\Z)=\Z  &  H^1(X_\epsilon,\Z)=\Z  \\[2mm]
H^3(X_\epsilon,\Z)=\Z  &  H^2(X_\epsilon,\Z)=\Z \oplus
\Coker(1-A_\epsilon) 
\end{array}
\end{equation} 
\endproof

\subsection{Action on $\Z^2$}\label{SZ2action}

We recall the following description of the action of 
$A_\epsilon$ on $\Lambda$, which will be useful in the
following, where we use twisted group $C^*$-algebras to
describe noncommutative tori.

\begin{lem}\label{actV1w}
In the basis $\{ 1, \theta \}$ of $\iota_1(L)\subset \R$, 
the action of the group $V=\epsilon^\Z$ 
is generated by the matrix 
\begin{equation}\label{varphie}
\varphi_\epsilon =\left(\begin{matrix} a & b \\ c & d
\end{matrix}\right) \in \SL_2(\Z)
\end{equation}
with $\epsilon = a + b \theta$ and $\epsilon\theta=c+d\theta$.
The conjugate elements $1/\theta$ and $1/\theta'$ are the fixed points of
$\varphi_\epsilon \in \SL_2(\Z)$ acting on $\P^1(\R)$ by fractional linear
transformations. 
\end{lem}

\proof As we have seen in \eqref{Uplusact}, the action 
of $V$ on $\Lambda$ is given by 
$$ A_\epsilon: (n+m\theta,n+m\theta') \mapsto (\epsilon (n+m\theta),\epsilon'
(n+m\theta')) $$
with $\epsilon\epsilon'=1$. In particular, for for $m=0$ 
and $n=1$ this gives $\epsilon \in \iota_1(L)$ and 
$\epsilon'\in \iota_2(L)$. Thus, we can write 
$\epsilon = a+b\theta$, for two integers $a,b\in\Z$. 
Similarly, the element $\epsilon \theta$ can be written 
in the form $\epsilon \theta =c+d\theta$.
Thus, the action 
of $A_\epsilon^k$ on $\Lambda$ can be described equivalently as
\begin{equation}\label{varphi}
 (n,m) \mapsto (n,m)\varphi_\epsilon^k , \ \ \ \text{ with } 
\varphi_\epsilon =\left(\begin{matrix} a & b \\ c & d 
\end{matrix}\right) \in \SL_2(\Z)
\end{equation}
The second statement follows immediately since 
$$\theta^{-1} =\epsilon/ (\epsilon\theta) = (a\theta^{-1}+b)/(c\theta^{-1}+d).$$
\endproof

We obtain in this way two corresponding identifications $S(\Lambda,V)=
\Z^2\rtimes_{\varphi_\epsilon}\Z$, by mapping $(\lambda,\epsilon^k)$
to either $(\lambda_1=n+m\theta,k)$ or $(\lambda_2=n+m\theta',-k)$.

\subsection{Solvmanifold and Hecke's lift of geodesics}\label{HeckeliftSec}

For $\Gamma=\SL_2(\Z)$ and $X_\Gamma=\Gamma\backslash \H$ the modular curve, let 
$\cU_\Gamma \to X_\Gamma$ denote the universal family of elliptic curves over the 
modular curve, where the fiber over $\tau\in X_\Gamma$ of $\cU_\Gamma$ is the isomorphism 
class of the elliptic curve $E_\tau=\C/(\Z+\tau\Z)$. 

Suppose given a lattice $L$ in a real quadratic field $\K$ and let $\{1,\theta\}$ be 
a basis for $L$, with $\theta'$ the Galois conjugate of $\theta$ under the Galois
involution of $\K$ over $\Q$. 

We denote by $\gamma_{\theta,\theta'}$ the infinite geodesic in the hyperbolic plane $\H$ 
with endpoints $\theta,\theta'\in \P^1(\R)$. This defines a closed geodesic in the
quotient $X_\Gamma$ of length $\log \epsilon$, for $\epsilon>1$ the generator of 
$V=U^+_L=\epsilon^\Z$. We denote the closed geodesic by $\bar\gamma_{\theta,\theta'}$.

Consider the restriction of the universal family $\cU_\Gamma$ to the closed geodesic
$\bar\gamma_{\theta,\theta'}$. Via the parameterization of the closed geodesic by
a circle $S^1$ of length $\log\epsilon$, we can consider the pullback to the 
parameterizing $S^1$ of $\cU_\Gamma$. We obtain in this way a real 3-dimensional manifold,
which we denote $\cU_{\theta,\theta'}$. This is topologically a 3-manifold that 
fibers over a circle, with fibers $T^2$. We consider it endowed with the metric that
is the product of the geodesic length and the flat metric on $T^2$.
We then obtain the following result.

\begin{lem}\label{3mfldUnivFam}
The solvmanifold $S(\Lambda,V)$ is isometrically 
equivalent to $\cU_{\theta,\theta'}$.
\end{lem}

\proof We recall the following procedure of Hecke to lift closed geodesics to the space
of lattices (\cf Manin \cite{Man}, \S 1.8.2).
Given a lattice $L$ in a real quadratic field $\K$, with $\ell \to \ell'$ the
Galois involution, one sets
\begin{equation}\label{Lambdat}
\Lambda_t(L):=\{ z\in \H\,|\, z=z(\ell,t)= \ell e^{t} + i\ell' e^{-t}\,\, \ell\in L \}.
\end{equation}
This defines, for all $t\in\R$ a lattice $\Lambda_t \subset \C$. The action of 
$V=\epsilon^\Z$ is of the form (\cf Lemma 1.8.3 of \cite{Man})
\begin{equation}\label{VactHecke}
z(\ell,t)\mapsto \epsilon\ell e^t+i \epsilon'\ell' e^{-t} =z(\ell,t+\log\epsilon).
\end{equation}
In particular (see again \cite{Man} Lemma 1.8.3), for $\{1,\theta\}$ a basis of $L$, 
the lattice $\Lambda_t(L)$ is generated by $\{ 1, \tau_t\}$ where $\tau_t$ runs over 
the geodesic $\ell_{\theta,\theta'}\subset \H$, for $t\in \R$. Thus, we can identify
the 3-manifold $\cU_{\theta,\theta'}$ with the fibration over a circle of length 
$\log\epsilon$, with fiber $E_{\tau_t}=\C/\Lambda_t(L)$.  

On the other hand, the 3-manifold $S(\Lambda,V)$ is a fibrations of tori 
over the circle
\begin{equation}\label{torifibration}
T^2 \to S(\Lambda,V) \to S^1,
\end{equation}
where the base $S^1$ is a circle of length $\log\epsilon$ and the 
fiber over $t\in S^1$ is given by the 2-torus
\begin{equation}\label{2torus}
T^2_t = \R^2/\Lambda_t,
\end{equation}
with $\Lambda_t=\Theta_t(\Lambda)$, for $\Theta_t(x,y)=(e^t x,e^{-t} y)$. 
This proves the result.
\endproof

\section{Actions on noncommutative tori with real multiplication}\label{Snctori}

The noncommutative torus $\cA_\theta$ of modulus 
$\theta \in \R\smallsetminus \Q$ is the noncommutative space
described, at the topological level, by the irrational rotation
$C^*$-algebra, that is, the universal $C^*$-algebra generated 
by two unitaries $U$, $V$ with the commutation relation $VU =e^{2\pi i
\theta} UV$. It has a smooth structure given by the smooth subalgebra
of series $\sum_{n,m} a_{n,m} U^n V^m$ with rapidly decaying
coefficients (\cf \cite{CoCR}). 

It is a well known result (\cite{CoCR}, \cite{Rie}) that the
algebras $\cA_{\theta_1}$ and $\cA_{\theta_2}$ are Morita equivalent
whenever there exists an element $g\in \SL_2(\Z)$ acting on $\R$ by
fractional linear transformations, such that $\theta_1=g\theta_2$. 
In the following we concentrate on the case where the irrational
number $\theta$ is a quadratic irrationality in a real quadratic field
$\K=\Q(\theta)$. These are the noncommutative tori with real
multiplication considered in \cite{Man}.
We let $L$ be the lattice in $\K$ with
$\iota_1(L)=\Z+\Z\theta$ and $\iota_2(L)=\Z+\Z\theta'$. As before,
we denote by $\Lambda$ the corresponding lattice in $\R^2$.

The $C^*$-algebra of the noncommutative torus $\cA_\theta$ described
above can be equivalently described as the crossed product
\begin{equation}\label{NCtorusS1}
\cA_\theta = C(S^1)\rtimes_\theta \Z, 
\end{equation}
where the action of $\Z$ on $S^1$ is by the irrational rotation by
$\exp (2\pi i \theta)$. Up to Morita equivalence, one can replace
$C(S^1)$ by the crossed product $C_0(\R)\rtimes \Z$, and one obtains a
Morita equivalent description of the noncommutative torus as
\begin{equation}\label{NCtorusR}
 C_0(\R)\rtimes_\theta \Z^2 = C_0(\R)\rtimes (\Z +
\Z \theta).
\end{equation}
In the case we are considering, of real quadratic fields, we can
regard the noncommutative torus with real multiplication associated to
a lattice $L\subset \K$ as described by the algebras 
\begin{equation}\label{NCtorusRrm}
\bT_{\Lambda,1}:= C_0(\R)\rtimes \iota_1(L) \ \ \ \ \ 
\bT_{\Lambda,2}:= C_0(\R)\rtimes \iota_2(L). 
\end{equation}
These algebras can be described as follows. They are $C^*$-algebras 
generated by elements of the form $f U_\lambda$, with $f\in C_0(\R)$ and 
$\lambda\in \Lambda$, with the product
$$ f U_\lambda\, h U_\eta = f U_{\lambda,i}(h) U_{\lambda+\eta}, \ \ \
\text{ where } \ \ 
 U_{\lambda,i}(h)(x)= h(x+\iota_i(\lambda)), \ \ \ i=1,2. $$

The group $V=\epsilon^\Z$ of units acts as
symmetries of the noncommutative tori $\bT_{\Lambda,i}$
as follows.

\begin{lem}\label{AutNCtori}
For $k\in\Z$ and $fU_\lambda \in \bT_{\Lambda,i}$, set
$\upsilon_\epsilon^k(f)(x) := f(\epsilon^k x)$ and
\begin{equation}\label{uepsilonk}
\upsilon_1^k (f U_\lambda) = \upsilon_\epsilon^k(f)
U_{A_\epsilon^k(\lambda)}, \ \ \ \text{ and } \ \ 
 \upsilon_2^k(f U_\lambda) = \upsilon_{\epsilon'}^k(f)
U_{A_\epsilon^k(\lambda)}.
\end{equation}
This defines actions $\upsilon_i : V \to
\Aut(\bT_{\Lambda,i})$.
\end{lem}

\proof The result follows directly from
$$ \upsilon_i^k(U_\lambda (h))(x) = U_{A_\epsilon^k(\lambda)}
(\upsilon_i^k(h))(x)= \left\{\begin{array}{ll} 
h(\epsilon^k(x+n+m\theta)) & i=1 \\
h(\epsilon^{-k}(x+n+m\theta')) & i=2 \end{array} \right. $$
which implies that
$$ \upsilon_\epsilon^k(f U_\lambda h U_\eta)=
\upsilon_i^k(f) \upsilon_i^k(U_\lambda(h)) 
U_{A_\epsilon^k(\lambda+\eta)} 
=\upsilon_i^k(f U_\lambda) \upsilon_i^k(h U_\eta). $$
\endproof

It is customary, in noncommutative geometry, to replace quotients
by crossed product algebras. In this case, the quotient
of the noncommutative tori $\bT_{\Lambda,i}$ by the action of
$V$ is described by the crossed product algebra
\begin{equation}\label{crossthetai}
\bT_{\Lambda,V,i} := \bT_{\Lambda,i} \rtimes_{\upsilon_i} V ,
\end{equation}
which we can view equivalently as the crossed product
\begin{equation}\label{crossthetai2}
\bT_{\Lambda,V,i} := C_0(\R)\rtimes_i  S(\Lambda,V),
\end{equation}
for the actions of $S(\Lambda,V)$ on $C_0(\R)$ of the form
\begin{equation}\label{actionZ2Z}
 U_{(\lambda,k)} f (x)= f(\epsilon^k (x+n+m\theta)) \ \ \text{ or } \ \ 
 U_{(\lambda,k)} f (x)= f(\epsilon^{-k} (x+n+m\theta')). 
\end{equation}

\section{Twisted group algebras and the magnetic Laplacian}\label{Stwist}

Another equivalent description of the algebra $\cA_\theta$ of the
noncommutative torus is as twisted group $C^*$-algebra. This played 
an important role in the context of the noncommutative geometry model 
of the integer quantum Hall effect (see \cite{Bel1}).

We recall briefly the definition and properties of twisted group
$C^*$-algebras, as this will be useful in the following. For a similar
overview and applications to the case of Fuchsian groups 
see \cite{MarMat3}. 

\subsection{Twisted group algebras}\label{Stwistalg}

Let $\Gamma$ be a finitely generated discrete group, and let $\sigma:
\Gamma\times \Gamma \to U(1)$ be a multiplier, that is, a 2-cocycle
satisfying the cocycle property
\begin{equation}\label{sigmacocycle}
\sigma(\gamma_1,\gamma_2)\sigma(\gamma_1\gamma_2,\gamma_3)=
\sigma(\gamma_1,\gamma_2\gamma_3)\sigma(\gamma_2,\gamma_3),
\end{equation}
with $\sigma(\gamma,1)=\sigma(1,\gamma)=1$. 

Consider then the Hilbert space $\ell^2(\Gamma)$ and the left/right
$\sigma$-regular representations of $\Gamma$ given by
\begin{equation}\label{LRsigmareps}
L^\sigma_\gamma f(\gamma')= f(\gamma^{-1}\gamma') \sigma(\gamma,
\gamma^{-1}\gamma'),  \ \ \  R^\sigma_\gamma f (\gamma')=
f(\gamma'\gamma) \sigma(\gamma',\gamma).
\end{equation}
They satisfy the relations
\begin{equation}\label{LRrepsprod}
L^\sigma_\gamma L^\sigma_{\gamma'}= \sigma(\gamma,\gamma')
L^\sigma_{\gamma\gamma'}, \ \ \  R^\sigma_\gamma R^\sigma_{\gamma'} =
\sigma(\gamma,\gamma') R^\sigma_{\gamma\gamma'}.
\end{equation}
Moreover the left $\sigma$-regular representation commutes with the
right $\bar\sigma$-regular representation, with $\bar\sigma$ the
conjugate multiplier. 
The algebra generated by the $R^\sigma_\gamma$ is the twisted group
ring $\C(\Gamma,\sigma)$. Its norm closure is the (reduced) twisted
group $C^*$-algebra $C^*_r(\Gamma,\sigma)$.

\subsection{The noncommutative tori as twisted group algebras}
\label{NCtorustwistSec}

One identifies the $C^*$-algebra $\cA_\theta$ of the noncommutative
torus with the reduced twisted
group $C^*$-algebra $C^*_r(\Z^2,\sigma)$ in the following way. Let
$\sigma$ be a cocycle of the form
\begin{equation}\label{sigmaZ2}
\sigma((n,m),(n',m')):= \exp(-2\pi i(\xi_1 nm' + \xi_2 mn')). 
\end{equation}
Then the operators
$U=R^\sigma_{(0,1)}$ and $V=R^\sigma_{(1,0)}$ acting by
$$ U f(n,m)=e^{-2\pi i \xi_2 n} \, f(n,m+1), \ \ \ 
V f(n,m)= e^{-2\pi i \xi_1 m}\, f(n+1,m) $$
that generate the algebra $C^*_r(\Z^2,\sigma)$
satisfy the commutation relation
$$ UV = e^{2\pi i \theta} VU, \ \ \ \text{ with } 
\theta =\xi_2-\xi_1. $$

Notice that different choices of $\xi_1,\xi_2$ with the same 
$\theta=\xi_2-\xi_1$ yield the same algebra $\cA_\theta$. 
This gives us the freedom to choose the $\xi_i$ according to the
following result. 

\begin{lem}\label{sigmaSL2inv}
A cocycle $\sigma$ of the form \eqref{sigmaZ2} 
has the property that
\begin{equation}\label{SL2inv}
\sigma((n,m),(n',m'))=\sigma((n,m)\varphi,(n',m')\varphi), \ \ \ 
\forall \varphi=\left(\begin{matrix} a&b\\c&d \end{matrix}\right) \in
\SL_2(\Z) 
\end{equation}
if and only if $\xi_2=-\xi_1$.
\end{lem}

\proof We see that $\sigma((n,m)\varphi,(n',m')\varphi)$ is of the form
$$ \exp(-2\pi i((\xi_1+\xi_2) (ab\, nn'+cd\,mm') +
(\xi_1 cb + \xi_2 ad) mn' + (\xi_1 ad +\xi_2 cb) nm')). $$
\endproof

Thus, in the following we will assume that $\xi_2 =\theta/2=-\xi_1$ in
the choice of the cocycle $\sigma$ of \eqref{sigmaZ2}. We can then write
$\sigma$ in the form
$$ \sigma_\theta ((n,m),(k,r))=\exp(\pi i \theta (nr-mk)) =
\exp(\pi i \theta (n,m)\wedge (k,r)), $$
where we use the notation
\begin{equation}\label{detwedge}
 (a,b)\wedge (c,d)= \det\left(\begin{matrix} a & b \\ c & d
\end{matrix}\right) .
\end{equation}

We then obtain the following identifications.

\begin{cor}\label{Lambdasigma}
The noncommutative tori $\bT_{\Lambda,i}$ are described by twisted
group $C^*$-algebras  
\begin{equation}\label{Lambdasigmai}
\begin{array}{rll}
\bT_{\Lambda,1}= & C^*(\Z^2, \sigma_\theta) & = C^*(\Lambda,
\sigma_{\theta(\theta'-\theta)^{-1}}), \\[2mm] 
\bT_{\Lambda,2}= & C^*(\Z^2, \sigma_{\theta'}) & = C^*(\Lambda,
\sigma_{\theta'(\theta'-\theta)^{-1}}). \end{array}
\end{equation}
\end{cor}

\proof The expression 
$$ \sigma_u (\lambda,\eta) =\exp(\pi i u\, \lambda\wedge \eta) $$
defines a cocycle on $\Lambda$. For $\lambda=(n+m\theta,n+m\theta')$
and $\eta=(k+r\theta,k+r\theta')$, a direct calculation shows that
$$ \sigma_\theta ((n,m),(k,r))= \sigma_u (\lambda,\eta), \ \ \ 
\text{ for } \  u=\theta (\theta' -\theta)^{-1}. $$
Thus, the generators 
$R^\sigma_{(n,m)}$ of $C^*(\Z^2, \sigma_\theta)$ with
$$ R^\sigma_{(n,m)} R^\sigma_{(k,r)} =\sigma_\theta ((n,m),(k,r))
R^\sigma_{(n,m)+(k,r)} $$
are identified with the generators $R^\sigma_\lambda$ of $C^*(\Lambda,
\sigma_{\theta(\theta'-\theta)^{-1}})$ with
$$ R^\sigma_\lambda R^\sigma_\eta = \sigma_{\theta(\theta'-\theta)^{-1}}
(\lambda,\eta) R^\sigma_{\lambda+\eta}. $$
The case of $\bT_{\Lambda,2}$ is analogous.
\endproof

\subsection{Twisted group algebra of $S(\Lambda,V)$}\label{StwistV}

We now show that the algebra $\cA_\theta \rtimes V$, 
which we introduced in the previous section to describe 
the quotient of the noncommutative torus
with real multiplication by the action of $V$, also admits a
description in terms of twisted group $C^*$-algebras, for
the group $S(\Lambda,V)$. First notice that the group 
$S(\Lambda,V)$ is amenable, so that the
maximal and reduced group $C^*$-algebras coincide, 
$C^*_{max}(S(\Lambda,V))\cong C^*_r(S(\Lambda,V))$,
so that we can simply write $C^*(S(\Lambda,V))$, and
$C^*(S(\Lambda,V),\tilde\sigma)$ for the twisted case.

\begin{lem}\label{tildesigmalem}
Let $\sigma$ be a multiplier on $\Z^2$ of the form \eqref{sigmaZ2},
with $\xi_2=\theta/2=-\xi_1$. Then the map $\tilde\sigma:
S(\Lambda,V)\times S(\Lambda,V)\to U(1)$ of the form
\begin{equation}\label{tildesigma}
\tilde\sigma((n,m,k),(n', m',k')):=
\sigma((n,m),(n',m')\varphi_\epsilon^k)  
\end{equation}
is a multilier for $S(\Lambda,V)$, identified with the group
$\Z^2\rtimes_{\varphi_\epsilon}\Z$.
\end{lem}

\proof The cocycle condition for $\sigma$ and the $\SL_2(\Z)$-invariance 
$\sigma((n,m)\varphi,(n',m')\varphi)=\sigma((n,m),(n'm'))$ imply that
$\tilde\sigma$ also satisfies the cocycle condition \eqref{sigmacocycle},
since we have
$$ \sigma((n_1,m_1),(n_2,m_2)
\varphi_\epsilon^{k_1})\sigma((n_1,m_1)+(n_2,m_2)\varphi_\epsilon^{k_1},
(n_3,m_3)\varphi_\epsilon^{k_1+k_2}) = $$
$$ \sigma((n_1,m_1),(n_2,m_2)\varphi_\epsilon^{k_1}
+(n_3,m_3)\varphi_\epsilon^{k_1+k_2})\sigma((n_2,m_2)\varphi_\epsilon^{k_1},
(n_3,m_3)\varphi_\epsilon^{k_1+k_2}). $$
\endproof

We then have the following result.

\begin{prop}\label{AthetaVCSsigma}
The algebras $\bT_{\Lambda,V,i}=\bT_{\Lambda,i}\rtimes_i V$ are
isomorphic to the algebras
\begin{equation}\label{LambdaVsigmai}
\begin{array}{rll}
\bT_{\Lambda,V,1}= & C^*(\Z^2, \sigma_\theta)\rtimes_{\upsilon_1} \Z 
& = C^*(\Z^2\rtimes_{\varphi_\epsilon}\Z,\tilde\sigma_\theta) \\[2mm]
= &  C^*(\Lambda, \sigma_{\theta(\theta'-\theta)^{-1}})
\rtimes_{\upsilon_1}
V & = C^*(S(\Lambda,V),\tilde\sigma_{\theta(\theta'-\theta)^{-1}}) , \\[2mm] 
\bT_{\Lambda,V,2}= & C^*(\Z^2, \sigma_{\theta'})\rtimes_{\upsilon_2}
\Z  & = C^*(\Z^2\rtimes_{\varphi_\epsilon}\Z,\tilde\sigma_{\theta'})
\\[2mm] 
= & C^*(\Lambda,\sigma_{\theta'(\theta'-\theta)^{-1}})\rtimes_{\upsilon_2}
V & = 
C^*(S(\Lambda,V),\tilde\sigma_{\theta'(\theta'-\theta)^{-1}}). \end{array}
\end{equation}
\end{prop}

\proof We just show explicitly one of the identifications. The others
follow similarly. The twisted group algebra 
$C^*(\Z^2\rtimes_{\varphi_\epsilon}\Z,\tilde\sigma_\theta)$
is generated by elements $R^{\tilde\sigma}_{(n,m,k)}$ satisfying 
$$ R^{\tilde\sigma}_{(n,m,k)} R^{\tilde\sigma}_{(n',m',k')} =
\tilde\sigma((n,m,k),(n',m',k'))
R^{\tilde\sigma}_{(n,m,k)(n',m',k')}. $$
The crossed product $C^*(\Z^2,\sigma_\theta)\rtimes_{\upsilon_1}\Z$ 
is generated by elements of the form $R^\sigma_{(n,m)}
\upsilon_\epsilon^k$. The map $R^{\tilde\sigma}_{(n,m,k)} 
\mapsto R^\sigma_{(n,m)} \upsilon_\epsilon^k$ gives an identification
of the generators, which also satisfies
$$ \begin{array}{rl}
R^\sigma_{(n,m)}\upsilon_\epsilon^k R^\sigma_{(n',m')}
\upsilon_\epsilon^{k'}= & R^\sigma_{(n,m)}
R^\sigma_{(n',m')\varphi_\epsilon^k}\, \upsilon_\epsilon^{k+k'}
\\[2mm] = & \sigma((n,m),(n',m')\varphi_\epsilon^k) \,
R^\sigma_{(n,m)+(n',m')\varphi_\epsilon^k}
\upsilon_\epsilon^{k+k'}. \end{array} $$
This gives an isomorphism $C^*(\Z^2,\sigma_\theta)\rtimes_{\upsilon_1}\Z
=C^*(\Z^2\rtimes_{\varphi_\epsilon}\Z,\tilde\sigma_\theta)$.
\endproof

\subsection{The magnetic Laplacian}\label{Smagnet}

Consider the general setting of a finitely generated 
discrete group $\Gamma$ acting on a
contractible space $\tilde X$ with compact quotient $X=\tilde
X/\Gamma$. Assume everything happens in the smooth category and we
think of $\tilde X$ as endowed with a metric that is invariant under
the action of $\Gamma$. Upon
choosing a base point $x_0\in \tilde X$, we can think of the discrete
set $\Gamma x_0$ as a crystal of charged ions and consider the
electron--ion interaction problem in $\tilde X$. This means that
electrons move in $\tilde X$ subject to a periodic potential. Under
resonable assumptions, one can make an {\em independent electron
approximation} and replace the $N$-particle Hamiltonian with the 
unbounded periodic electric potential of the ion crystal with a single
electron Hamiltonian in an effective periodic potential given by a
bounded function (see \cite{MarMat3} for a brief overview).  

The Hamiltonian is then of the
form $\Delta + V$, where the $\Delta$ is the Laplacian on
$\tilde X$. We think of it as an unbounded operator on $L^2(\tilde
X)$. The Hamiltonian is invariant under translations by
$\gamma\in\Gamma$, that is, $T_\gamma \Delta = \Delta T_\gamma$ and by
construction $V$ is also invariant. Here the $T_\gamma$ are the
operators $T_\gamma f(x)=f(x\gamma)$ on $L^2(\tilde X)$. 

One can consider on $\tilde X$ a magnetic field. This is specified by
a closed 2-form $\omega$ which satisfies $\gamma^*\omega =\omega$. Since
$\tilde X$ is contractible, there is a global magnetic potential
$\omega =d\chi$. The corresponding hermitian connection $\nabla = d -i
\chi$ satisfies $\nabla^2 = i \omega$. 
The invariance of $\omega$ implies
$d(\chi-\gamma^*\chi)=0$, so that $\chi-\gamma^*\chi = d\phi_\gamma$,
where the function $\phi_\gamma(x)=\int_{x_0}^x \chi - \gamma^*\chi$ 
has the properties that
$$ 
\phi_\gamma(x) -\phi_{\gamma'}(\gamma x) -\phi_{\gamma\gamma'}(x)
$$
is independent of $x\in \tilde X$ and $\phi_\gamma(x_0)=0$ so that
\begin{equation}\label{cocyclephigamma}
 \sigma(\gamma,\gamma')=\exp(-i\phi_\gamma(\gamma'x_0))
\end{equation}
defines a multiplier $\sigma:\Gamma\times\Gamma\to U(1)$. 
The Laplacian $\Delta$ is naturally replaced, in the presence of a
magnetic field, by the magnetic Laplacian $\Delta^\chi
=\nabla^*\nabla= (d-i\chi)^*(d-i\chi)$. This is no longer invariant
under translations $T_\gamma$, but is invariant under the magnetic
translations
\begin{equation}\label{invTphigamma}
 T_\gamma^\phi \Delta^\chi = \Delta^\chi T^\phi_\gamma
\end{equation}
where $T_\gamma^\phi f(x)= \exp(i\phi_\gamma(x)) f(\gamma^{-1}x)$.
Similarly, in the independent electron approximation, the effective 
potential $V$ is also invariant under the magnetic translations.
  
The magnetic translations satisfy the relations of the twisted group
algebra $C^*_r(\Gamma,\bar\sigma)$
$$ T^\phi_\gamma T^\phi_{\gamma'} =\bar\sigma(\gamma,\gamma')
T^\phi_{\gamma\gamma'} , $$
for $\sigma$ as in \eqref{cocyclephigamma} and $\bar\sigma$ the conjugate. 
(We refer the reader to \cite{MarMat1}, \cite{MarMat3} for a brief overview 
of these well known facts.)

\subsection{Discretized electron--ion interaction and Harper operators}\label{Sharper}

It is usually convenient to discretize the electron--ion 
interaction problem. This means replacing the continuum model
with Hilbert space $L^2(\tilde X)$ by a discrete model on the
Hilbert space $\ell^2(\Gamma)$.
In the case without magnetic field, this is done
by replacing the Laplacian $\Delta$ by its discretized version 
$\Delta_{discr} = r-\cR$, where $r$ is the cardinality 
of a symmetric set of generators
for $\Gamma$ and $\cR$ is the random walk operator
\begin{equation}\label{randomwalk}
\cR = \sum_{i=1}^r R_{\gamma_i} \ \ \ \text{ with } \ \ \ 
R_{\gamma_i} f(\gamma)= f(\gamma\gamma_i) 
\end{equation}
for $f\in \ell^2(\Gamma)$. As in the continuum model the discretized
Laplacian commutes with translations by elements $\gamma\in
\Gamma$. The effective potential is then taken to be an element in the
group ring $\C[\Gamma]$. 

In the presence of a magnetic field, one can
still obtain a good discretized version of the electron--ion
interaction problem as in \cite{Sunada}. The random walk operator
of \eqref{randomwalk} is then replaced by the Harper operator
\begin{equation}\label{Harper}
\cH_\sigma = \sum_{i=1}^r R^\sigma_{\gamma_i},
\end{equation}
with $R^\sigma_{\gamma_i}\in \C(\Gamma,\sigma)$ the elements of the 
right $\sigma$-regular representation, with $\sigma$ the cocycle of
\eqref{cocyclephigamma}. The discretized version of the magnetic
Laplacian is then given by the operator
\begin{equation}\label{magnLapD}
\Delta^\chi_{discr} = r -\cH_\sigma,
\end{equation} 
which commutes with the magnetic translations 
$L_\gamma^{\bar\sigma}$. Similarly the effective potential is taken to
be an element $V\in \C(\Gamma,\sigma)$, which then also commutes with
the magnetic translations $L^{\bar\sigma}_\gamma$.

\subsection{Harper operators for noncommutative tori and for
$S(\Lambda,V)$}\label{SharperV}

In the case of the noncommutative torus, viewed as the twisted group
$C^*$-algebra $C^* (\Z^2,\sigma)$, the Harper operator is of the form
\begin{equation}\label{HarperAtheta}
\cH_\sigma = U + U^* + V + V^*,
\end{equation}
where $U$ and $V$ are the generators of $\cA_\theta$.

The spectral theory of the Harper operator $\cH_\sigma$ 
of \eqref{HarperAtheta} was widely studied. 
In particular, it was shown in \cite{hofs}
that the spectrum exhibits a remarkable fractal structure (the
Hofstadter butterfly) which appears to have infinitely many gaps
(Cantor like spectrum) for irrational $\theta$ and finitely many gaps
(band spectrum) for rational $\theta$. 
The precise gap structure of the spectrum of
Harper operators, as a function of the magnetic flux (that is $\theta$
in the noncommutative torus case), is a problem still under active
investigation. As we see more in detail in the following, in the
specific case of interest here, the gap labelling problem
for the spectrum of the Harper operator is closely related to the 
computation of the range of the trace on the $K$-theory of the twisted
group $C^*$-algebra. 

In the following, we will be interested in the case of the group
$S(\Lambda,V)$. In this case, after identifying it with
$\Z^2\rtimes_{\varphi_\epsilon}\Z$, the Harper operator is of the form
\begin{equation}\label{harperSLV}
\cH_{\tilde\sigma} = U + U^* + V + V^* + W + W^*,
\end{equation}
where $U=R^{\tilde\sigma}_{(0,1,0)}$, $V=R^{\tilde\sigma}_{(1,0,0)}$
and $W=R^{\tilde\sigma}_{(0,0,1)}$.

\subsection{Spectral theory and $K$-theory}\label{SspectralK}

We recall here briefly the relation between spectral theory of Harper
operators and $K$-theory of twisted group $C^*$-algebras (\cf
\cite{Bel3}, \cite{MarMat3} \S 3). We then
proceed in the following section to analyze the specific case of 
$C^* (S(\Lambda,V))$.

As we have seen, the twisted group $C^*$-algebra
$C_r^*(\Gamma,\sigma)$ is the norm closure of the twisted group ring
$\C(\Gamma,\sigma)$ in the right $\sigma$-regular representation on
$\ell^2(\Gamma)$, that is, the $C^*$-algebra generated by the magnetic
translations $R^\sigma_\gamma$. If we take the weak closure of
$\C(\Gamma,\sigma)$, we obtain the twisted group von Neumann 
algebra $\cU(\Gamma,\sigma)$. Suppose given an operator $\cH_{\sigma,
V}=\cH_\sigma +V$, with $\cH_\sigma$ the Harper operator described
above and $V$ an effective potential in $\C(\Gamma,\sigma)$. We have
by construction $\cH_{\sigma, V}\in \C(\Gamma,\sigma)\subset
C_r^*(\Gamma,\sigma)\subset \cU(\Gamma,\sigma)$, hence the spectral
projections of $\cH_{\sigma, V}$,
\begin{equation}\label{Specproj}
P_E={\bf 1}_{(-\infty,E]}(\cH_{\sigma, V})
\end{equation}
are in the von Neumann algebra, $P_E \in \cU(\Gamma,\sigma)$. In
particular, if the energy level $E$ is not in the spectrum of
$\cH_{\sigma,V}$, then the corresponding spectral projection $P_E$ is
actually in the $C^*$-algebra $C_r^*(\Gamma,\sigma)$. 

This implies that the question of counting gaps inthe spectrum of
$\cH_{\sigma,V}$ can be reformulated as a problem of counting
projections in the $C^*$-algebra $C_r^*(\Gamma,\sigma)$, modulo the
Murray-von Neumann equivalence relation, $P\sim Q$ if there exists
$V\in C_r^*(\Gamma,\sigma)\otimes \cK$ with $P= V^*V$
and $Q= VV^*$. Equivalent spectral projections correspond to a same gap in
the spectrum. The group $K_0(C_r^*(\Gamma,\sigma)))$ is the
Grothendieck group of the resulting abelian semi-group (with the
operation of direct sum). Thus, the gap counting problem is restated
as a problem involving $K$-theory of $C^*$-algebras. More precisely, 
there is a faithful canonical finite trace
$$ \tau(a) = \langle  a\delta_1, \delta_1\rangle_{\ell^2(\Gamma)}, $$
on $C_r^*(\Gamma,\sigma)$, with $\delta_\gamma$ the canonical basis of
$\ell^2(\Gamma)$. This extends to 
$$ \tr=\tau\otimes\Tr: \{ P\in C_r^*(\Gamma,\sigma)\otimes\cK)\,|\,
P^*=P,\ P^2=P \} \to \R, $$
with $\Tr$ the standard trace on bounded operators and induces
\begin{equation}\label{trK0}
[\tr]: K_0(C_r^*(\Gamma,\sigma)))\to \R.
\end{equation}
One can obtain an estimate of the number of equivalence classes of
projections by a direct computation of the range of the trace on
$K_0(C_r^*(\Gamma,\sigma)))$, using
\begin{equation}\label{trgapcount}
\tr(\{ P\in C_r^*(\Gamma,\sigma)\otimes\cK)\,|\,
P^*=P,\ P^2=P \})=
[\tr](K_0(C_r^*(\Gamma,\sigma)))\cap [0,1].
\end{equation}

\section{Homotopy quotient and the Baum--Connes conjecture}\label{SconjBC}

As we show in this section, the computation of the range of the trace
on $K$-theory is closely related to the use of the 3-manifold
$X_\epsilon$ as a commutative model up to homotopy of the
noncommutative space $\bT_{\Lambda,V,i}$. 

The main idea of the Baum--Connes conjecture is precisely the fact 
that noncommutative spaces originating from ``bad quotients''
have good homotopy quotients that can be used to compute geometrically
invariants such as the analytic K-theory. 

The group $S(\Lambda,V)$ we are considering here is a particular case 
of a class of groups of the form $\Z^2\rtimes_\varphi \Z$, for 
some $\varphi\in \SL_2(\Z)$. The corresponding (twisted) group 
$C^*$-algebras and their $K$-theory were analyzed in \cite{PaRae}. 
We wish to stress here the relation between the noncommutative space 
and its model $X_\epsilon$ and the role of the latter in the index
computations. 

\subsection{$K$-theory of $C^* (S(\Lambda,V),\tilde
\sigma)$}

We now compute explicitly the $K$-theory of the twisted group
$C^*$-algebra of $S(\Lambda,V)$. This can be done using the
Pimsner--Voiculescu six terms exact sequence.

\begin{lem}\label{KtheorySLVsigma}
The $K$-theory groups of $C^* (S(\Lambda,V),\tilde\sigma))$ are of
the form
\begin{equation}\label{K01groups}
\begin{array}{l}
K_0(C^* (S(\Lambda,V),\tilde\sigma))\cong \Lambda  \\[2mm]
K_1(C^* (S(\Lambda,V),\tilde\sigma)) \cong \Lambda\oplus
\Lambda/(1-A_\epsilon)\Lambda.
\end{array} 
\end{equation}
\end{lem}

\proof By Proposition \ref{AthetaVCSsigma}, we can identify
$C^* (S(\Lambda,V),\tilde\sigma)$ with the crossed product
$C^* (\Lambda,\sigma)\rtimes V$. Thus, we can apply the
Pimsner--Voiculescu six terms exact sequence for the actions of
$V\cong \Z$. We have
\begin{center}
\diagram
& & & K_0(\cA)\rto^{1-\alpha_*}& K_0(\cA)\rto & 
K_0(\cA\rtimes\Z) \dto^\partial \\ 
& & & K_1(\cA\rtimes\Z)\uto_\partial & 
K_1(\cA)\lto & K_1(\cA) \lto_{1-\beta_*} 
\enddiagram
\end{center}
where $\cA=C^* (\Lambda,\sigma)$ and $\alpha_*$ and $\beta_*$ denote,
the action on $K_0(\cA)$ and $K_1(\cA)$, respectively, induced by the
generator $A_\epsilon$ of the $\Z$-action on $\cA$. We can identify 
$K_0(\cA) = \Lambda = K_1(\cA)$. We then have $1-\alpha_* =0$
and $1-\beta_*=1-A_\epsilon$, so that we obtain
$\Ker(1-\beta_*)=\Ker(1-A_\epsilon)=0$ and
$\Coker(1-\beta_*)=\Lambda/(1-A_\epsilon)\Lambda$. 
\endproof

We find in this way an abstract isomorphism of abelian groups 
\begin{equation}\label{HiKi}
\begin{array}{lll}
K_0(C^* (S(\Lambda,V)),\tilde\sigma) & \cong 
H^{odd} (X_\epsilon,\Z) & \cong \Z^2  \\[2mm]
K_1(C^* (S(\Lambda,V)),\tilde\sigma) & \cong 
H^{ev} (X_\epsilon,\Z) & \cong \Z^2 \oplus
\Coker(1-A_\epsilon).
\end{array}
\end{equation}
This identification can be justified more naturally in terms of
the Baum--Connes conjecture, as we discuss in the following.

\subsection{$K$-theory and the twist}\label{StwistK}

The following result shows that the
presence of the twisting by $\tilde\sigma$ has no effect on 
the $K$-theory.

\begin{lem}\label{Kisom}
There is an isomorphism
\begin{equation}\label{KsigmaKiso}
K_i(C^* (S(\Lambda,V),\tilde\sigma))\cong K_i(C^* (S(\Lambda,V))
\end{equation}
between the $K$-theory of the twisted group $C^*$-algebra
$C^* (S(\Lambda,V),\tilde\sigma)$ and the $K$-theory of
the untwisted $C^* (S(\Lambda,V))$.
\end{lem}

\proof The argument is similar to that used in \cite{MarMat1},
\cite{MarMat2} and Corollary 2.2 of
\cite{ELPW}. 
The cocycle $\tilde\sigma$ is real in the sense of Definition 1.12 of
\cite{ELPW}, being of the form \eqref{tildesigma}, with $\sigma$ of
the exponential form $\sigma((n,m),(n',m'))= \exp(-\pi i\theta(mn'-
nm'))$. Thus, as in Corollary 1.13 of \cite{ELPW}, the identification
\eqref{KsigmaKiso} follows using a homotopy $\exp(-t\pi i \theta(mn'-
nm'))$, with $t\in [0,1]$.
\endproof

Notice in fact that for groups of the form $\Gamma=\Z^2\rtimes_\varphi \Z$,
with $\varphi\in \SL_2(\Z)$ all cocycles $\sigma:\Gamma\times\Gamma\to
U(1)$ are real in the above sense. This was observed already in
\cite{PaRae}.

\begin{lem}\label{realsigma}
Let $\sigma\in Z^2(\Gamma,U(1))$ be a cocycle. Then $\sigma$ is
cohomologous to a real cocycle, that is, to an element of
$Z^2(\Gamma,U(1))$ that is of the form $\exp(2\pi i \zeta)$ for
$\zeta \in Z^2(\Gamma,\R)$.
\end{lem}

\proof
We can see it easily as in \S 2.2 of \cite{MarMat1}, by
considering the exact sequence of coefficient groups
$$ 1 \to \Z \stackrel{\iota}{\to}\R 
\stackrel{\exp(2\pi i \cdot)}{\to} U(1) \to 1 $$
and the long exact cohomology sequence
\begin{equation}\label{cohseq} 
 \cdots \to H^2(\Gamma,\Z)\to H^2(\Gamma,\R)\stackrel{\exp(2\pi i
\cdot)_*}{\to} H^2(\Gamma,U(1))\stackrel{\delta}{\to}
H^3(\Gamma,\Z)\stackrel{\iota_*}{\to} H^3(\Gamma,\R) \to \cdots 
\end{equation}
Since in our case, for $\Gamma=S(\Lambda,V)$, we have $E\Gamma
=\R^2\rtimes \R$ and $B\Gamma=X_\epsilon$, we see that
$$ H^2(\Gamma,\Z)=H^2(X_\epsilon,\Z)=H_1(X_\epsilon,\Z)=\Z \oplus
\Lambda/(1-A)\Lambda $$
and
$$ H^3(\Gamma,\Z)= H^3(X_\epsilon) =\Z. $$
We then see that, in the sequence \eqref{cohseq} the map $\iota_*$ is
injective so that $\delta =0$. Thus, all elements in
$H^2(\Gamma,U(1))$ come from $H^2(\Gamma,\R)$ via the exponential
map. 
\endproof

In fact, we do not need this general fact, as the cocycle we are
using is already constructed in the desired exponential form,
but we stated it here for completeness.

\subsection{Thom isomorphism, homotopy quotients and Baum--Connes}\label{Sthomiso}

It is known that the group $S(\Lambda,V)$ satisfies the Baum--Connes
conjecture (with coefficients). In fact the group $\SL_2(\Z)$ is 
known to satisfy the Baum--Connes conjecture with coefficients, 
hence by \cite{ChaEch} so does the group
$\Z^2\rtimes_{\varphi_\epsilon}\Z$ with $\varphi_\epsilon\in
\SL_2(\Z)$.

This means that the Kasparov assembly map is an isomorphism, 
hence the $K$-theory of the $C^*$-algebra
$C^*(S(\Lambda,V))$ can be computed in terms of the geometric
$K$-theory of the homotopy quotient $\underline{B}\Gamma$, the
classifying space for proper actions (\cf \cite{BC2}). 
This relates directly the analytic $K$-theory of the $C^*$-algebra
to the topological $K$-theory of the 3-manifold $X_\epsilon$

\begin{lem}\label{mumap}
The Kasparov assembly map for $C^*(S(\Lambda,V))$ gives an isomorphism
\begin{equation}\label{muiso}
\begin{array}{l}
\mu: K^1(X_\epsilon) \stackrel{\cong}{\to}
K_0(C^*(S(\Lambda,V))) \\[2mm]
\mu: K^0(X_\epsilon) \stackrel{\cong}{\to}
K_1(C^*(S(\Lambda,V))).
\end{array}
\end{equation}
\end{lem} 

\proof In our case the space $\underline{E}\Gamma$ is the solvable Lie group 
$S(\R^2,\R)=\R^2\rtimes \R$ and the homotopy
quotient $\Gamma \backslash \underline{B}\Gamma$ is the 3-manifold 
$X_\epsilon = S(\Lambda,V)\backslash S(\R^2,\R)$. This can be
identified with the mapping torus $$ X_\epsilon = T^2\times
[0,1]/((x,y),0)\sim (A_\epsilon (x,y), 1).$$
For a mapping torus, the Thom isomorphism \cite{Co-thom}  
gives the identification 
\begin{equation}\label{maptorus}
K_{i+1}(C(X_\epsilon))= K_i(C(T^2)\rtimes_{A_\epsilon}\Z).
\end{equation}
Moreover, the $C^*$-algebra $C(T^2)\rtimes_{A_\epsilon}\Z$
is identified with $C^*(\Lambda)\rtimes V$ by Fourier transform,
which identifies $C(T^2)=C^*(\Lambda)$ for $T^2=\R^2/\Lambda$.
The algebra $C^*(\Lambda)\rtimes V$ is then 
isomorphic to $C^*(S(\Lambda,V))$, by the same argument
of Proposition \ref{AthetaVCSsigma} in the untwisted case.
\endproof

\section{Twisted index theorem, $K$-theory, and the range of the
trace}\label{SrangeTr}

As we have seen, the 2-cocycle $\tilde\sigma$ on $\Gamma=S(\Lambda,V)$, is of
the form $\tilde\sigma = \exp(2\pi i \zeta)$, with $\zeta\in
H^2(\Gamma,\R)$. Upon identifying
$H^2(\Gamma,\R)=H^2(B\Gamma,\R)=H^2(X_\epsilon,\R)$, we can identify
the 2-cocycle $\zeta$ with a closed 2-form $\omega_\epsilon$ on the 3-manifold
$X_\epsilon$. We denote by $\tilde\omega_\epsilon$ its pullback to the
universal covering $\tilde X_\epsilon=S(\R^2,\R)$. This is a
$\Gamma$-invariant 2-form, $\gamma^*\omega_\epsilon =\omega_\epsilon$, 
which we previously interpreted as a magnetic field. 

\begin{lem}\label{lemzetaomega}
The real 2-cocycle $\zeta\in H^2(\Gamma,\R)$ with $\tilde\sigma
=\exp(2\pi i \zeta)$ is given by
\begin{equation}\label{zetaomega}
\zeta((\lambda,k),(\eta,r))=\frac{1}{4\pi i} \int_{\cR} \omega,
\end{equation}
where $\omega$ is the closed 2-form on $T^2=\R^2/\Lambda$ associated to the
cocycle $\sigma$ on $\Lambda$, with magnetic flux $\int_{T^2}\omega=2\pi
i\theta(\theta'-\theta)^{-1}$, and 
$\cR\subset \R^2$ is the oriented parallelogram with vertices
\begin{equation}\label{parallR}
\{ 0,A_\epsilon^k(\eta),\lambda,\lambda+A_\epsilon^k(\eta) \}.
\end{equation}
\end{lem}

\proof On $\tilde X_\epsilon$
the form $\tilde\omega_\epsilon$ is exact, hence we have a global magnetic
potential $\chi_\epsilon$ with
$\tilde\omega_\epsilon=d\chi_\epsilon$ and
$d(\chi_\epsilon-\gamma^*\chi_\epsilon)=0$, 
or $\chi_\epsilon -\gamma^*\chi_\epsilon =d\phi_\gamma$, as before, where the
$\phi_\gamma$ recovers the cocycle $\tilde\sigma$ by the formula
$$ \tilde\sigma(\gamma,\gamma') =\exp(-\phi_\gamma(\gamma'
x_0))=\exp(\int_{x_0}^{\gamma' x_0} \gamma^*\chi_\epsilon -\chi_\epsilon). $$

We know from Lemma \ref{tildesigmalem} that the cocycle $\tilde\sigma$
has the form $\tilde\sigma_u((\lambda,k), (\eta,r)) = \sigma_u(\lambda,
A_\epsilon^k (\eta))$, for $u=\theta(\theta'-\theta)^{-1}$,
so that we have
$$\zeta((\lambda,k),(\eta,r))= \frac{\theta(\theta'-\theta)^{-1}}{2}\, 
A_\epsilon^k (\eta) \wedge \lambda,$$
that is, $\zeta((\lambda,k),(\eta,r))
=\xi(\lambda,A_\epsilon^k(\eta))$, where 
\begin{equation}\label{zetaxi}
\xi(\lambda,A_\epsilon^k(\eta))
=\frac{1}{2\pi i} \int_0^{A_\epsilon^k(\eta)}
U_{\lambda}^*\chi -\chi.
\end{equation}
Here $\chi$ is the magnetic potential on $\R^2$ 
associated to the closed 2-form $\omega$ with
$$ \int_{T^2} \omega = 2\pi i \theta (\theta'-\theta)^{-1}. $$

Let then $\cR$ denote the oriented parallelogram in $\R^2$ with
vertices as in \eqref{parallR}.
We have 
\begin{equation}\label{intRomega}
\frac{1}{2\pi i}\int_{\cR} \omega =\frac{1}{2\pi
i}\int_0^{A_\epsilon^k(\eta)} (U_{\lambda}^* \chi -\chi) -
\frac{1}{2\pi i}\int_0^{\lambda} (U_{A_\epsilon^k(\eta)}^*
\chi -\chi).
\end{equation}
Using the fact that $\xi(\eta,\lambda)=-\xi(\lambda,\eta)$, 
this gives
$$ \frac{1}{2\pi i}\int_{\cR} \omega =2 \zeta((\lambda,k),(\eta,r)). $$
\endproof

\subsection{Spectral flow and odd Chern character}\label{Schern}

An element of $K_1(C(X_\epsilon))$ can be viewed as the class $[g]$ of
$g\in U_N(C(X_\epsilon))$, which we can see as a differentiable map
$g: X_\epsilon \to \GL_N(\C)$. We proceed as in \cite{Getz} and we
consider the associated 1-form 
\begin{equation}\label{betaform}
\beta(g)= g^{-1}dg \in \Omega^1(X_\epsilon, gl_N(\C)). 
\end{equation}
The corresponding family of connections
$\nabla_u = d + u\beta(g)$ on the trivial bundle $X_\epsilon\times
\C^N$ determines a closed Chern--Simons form
\begin{equation}\label{csform}
Ch(g):=cs(d,d+\beta(g))= \int_0^1 \Tr(\frac{d}{du}(\nabla_u) e^{\nabla_u^2})
du ,
\end{equation}
which gives the odd Chern character $Ch(g)$.
As shown in \cite{Getz}, this has an expression as an odd differential
form 
\begin{equation}\label{Chgbeta}
Ch(g)=\sum_{k=0}^\infty (-1)^k \frac{k!}{(2k+1)!} \Tr(\beta(g)^{2k+1}).
\end{equation}

One then has, see \cite{Getz}, that the pairing
\begin{equation}\label{Dgpair}
\langle D, [g]\rangle = SF(D,g^{-1}D g)
\end{equation}
of an odd Fredholm module $(\cH,D)$ with $[g]\in K_1$ is given by the
spectral flow along $D_u=(1-u) D +u g^{-1}D g$. In the case where
$D=\dirac$ is the Dirac operator of a compact spin manifold, this is computed
by the Atiyah--Patodi--Singer index formula \cite{APS}. In our case, 
this gives
\begin{equation}\label{indexDg}
SF(\dirac ,g^{-1}\dirac g)=-\frac{1}{(2\pi i)^2} \int_{X_\epsilon} \hat
A(X_\epsilon) Ch(g).
\end{equation}

\subsection{Twisted index theorem}\label{StwistInd}

We need the twisted version of \eqref{indexDg} above. 
Let $\dirac=\dirac_{X_\epsilon}$ be the Dirac operator 
on $X_\epsilon$, and let $\tilde\dirac$ be its lift to 
the universal cover $\tilde X_\epsilon = S(\R^2,\R)$. 
We then consider the twisting 
$\tilde\dirac \otimes \nabla$ of
the operator $\tilde\dirac$ by the hermitian connection $\nabla =d + i
\eta_\epsilon$ on the trivial line bundle on $\tilde X_\epsilon$, with
$\eta_\epsilon$ the 1-form giving the magnetic potential
$d\eta_\epsilon =\omega_\epsilon$ on $\tilde X_\epsilon$. 

While the operator $\tilde\dirac$ is $\Gamma$-invariant, 
with $\Gamma=S(\Lambda,V)$,
the twisted operator $\tilde\dirac \otimes \nabla$ is only invariant
under the projective action $(\Gamma,\tilde\sigma)$ of the magnetic
translations $R^{\tilde\sigma}_\gamma$. 

Consider then the 1-parameter family of operators $D_u=\tilde\dirac_u
\otimes \nabla$, where $\tilde\dirac_u=(1-u)\tilde\dirac + u g^{-1}
\tilde \dirac \, g$, for $[g]\in K^1(X_\epsilon)$ and
the associated operator $\cD_g =\frac{\partial}{\partial u} +  D_u$ on 
$\tilde X_\epsilon\times [0,1]$, which we can extend to $\tilde
X_\epsilon \times \R$ (\cf p.95 of \cite{APS}).

\begin{thm}\label{twistindthm}
The range of the trace on $K_0(C^*(S(\Lambda,V),\tilde\sigma))$ 
is given by
\begin{equation}\label{indrangetr}
[\tr](\mu_{\tilde\sigma}[g])= \frac{-1}{(2\pi
i)^2} \int_{X_\epsilon} \hat A \, e^{\omega_\epsilon} \, Ch(g), 
\end{equation}
where $\mu_{\tilde\sigma}:K^1(X_\epsilon)\to
K_0(C^*(S(\Lambda,V),\tilde\sigma))$ is the (twisted) Kasparov isomorphism,
$[g]\in K^1(X_\epsilon)$, and $\omega_\epsilon$ is the closed 2-form on
$X_\epsilon$ associated to the cocycle $\tilde\sigma$. 
\end{thm}

\proof
We let $P^\pm$ be the projections on the $L^2$-kernel of $\cD_g
\cD_g^*$ and $\cD_g^*\cD_g$, respectively, namely
$$ \cD_g P^+=0 \ \ \ \ \cD_g^* P^- =0 . $$

The $P^\pm$ have smooth kernels $P^\pm(x,y)$ and the
$(\Gamma,\tilde\sigma)$-invariance of $\tilde\dirac\otimes\nabla$
implies that
$$ e^{-i\phi_\gamma(x)} P^\pm(\gamma x,\gamma y) e^{i\phi_\gamma(y)}
=P^\pm(x,y), $$
which implies that $P^\pm(x,x)$ is $\Gamma$-invariant, for
$\Gamma=S(\Lambda,V)$.

We proceed as in \cite{At} and consider the von Neumann trace
$$ \tr (P^\pm) = \int_{X_\epsilon\times S^1} tr P^\pm((x,t),(x,t))\,
dx\, dt, $$
where $tr\,P^\pm(x,x)$ is the pointwise trace. 
The $L^2$-index of $\cD_g$ is given by 
\begin{equation}\label{L2ind}
\Ind_{L^2} (\cD_g)=\tr(P^+)-\tr(P^-).
\end{equation}

We define $\bar P^\pm$ by the smooth kernels
\begin{equation}\label{barPpm}
\bar P^\pm(x,y) = \int_{S^1} tr P^\pm ((x,t),(y,t))\, dt .
\end{equation}
These satisfy $\tr(\bar P^\pm)=\tr(P^\pm)$ by
$$ \int_{X_\epsilon\times S^1} tr P^\pm((x,t),(x,t))\,
dx\, dt =\int_{X_\epsilon} tr \bar P^\pm(x,x) dx. $$ 

The projections $\bar P^\pm$ are in the von Neumann algebra
$\cU(\Gamma,\tilde\sigma)$. After adding a compact perturbation 
in $C^*_r(\Gamma,\tilde\sigma)$ one obtains a well defined 
index (\cf \cite{MiFo}, \cite{Ros}), 
\begin{equation}\label{indGs}
\Ind_{(\Gamma,\tilde\sigma)}(\cD_g) =[\bar P^+]-[\bar P^-]\in
K_0(C^*_r(\Gamma, \tilde\sigma)).
\end{equation}
The (twisted) Kasparov map $\mu: K^1(X_\epsilon)\to
K_0(C^*(\Gamma,\tilde\sigma))$ is given by
\begin{equation}\label{KaspK1}
\mu_{\tilde\sigma}[g] = \Ind_{(\Gamma,\tilde\sigma)}(\cD_g).
\end{equation}

We obtain in this way that
$$ \Ind_{L^2}(\cD_g)=\tr(\bar P^+)-\tr(\bar P^-)=\tr
(\Ind_{(\Gamma,\tilde\sigma)}(\cD_g)). $$

Consider the heat kernel $e^{-t \cD^2}$, where 
$$ \cD = \left(\begin{matrix} 0 & \cD_g^* \\ \cD_g & 0
\end{matrix}\right) \ \ \ \text{ with } \ \ \
 \cD^2 =  \left(\begin{matrix} \cD_g^* \cD_g & 0 \\ 0 & \cD_g \cD_g^*
\end{matrix}\right). $$
We have
$$ \lim_{t\to \infty} \tr_s(e^{-t \cD^2})= \tr(P^+)-\tr(P^-) $$
and
$$ \frac{\partial}{\partial t} \tr_s(e^{-t \cD^2}) = -\tr_s(\cD^2
e^{-t \cD^2}) =\tr_s([\cD e^{-t \cD^2},\cD])=0. $$
Thus 
$$ \tr(P^+)-\tr(P^-) = \lim_{t\to \infty} \tr_s(e^{-t \cD^2})=
\lim_{t\to 0} \tr_s(e^{-t \cD^2}) $$
$$  = \frac{-1}{(2\pi
i)^2}\int_{X_\epsilon\times S^1} \hat A\,
Ch(\nabla_u) =\frac{-1}{(2\pi
i)^2}\int_{X_\epsilon} \hat A\, e^{\omega_\epsilon} \, Ch(g), $$
where $Ch(\nabla_u)=tr(\beta e^{(d+u\beta)^2})$ for $\beta=g^{-1}dg$,
with $\int_{S^1} Ch(\nabla_u)= Ch(g)$.
\endproof

\subsection{Range of the trace}\label{Srange}

Using the twisted index theorem we can then compute explicitly the range of
the trace on $K_0(C^*(S(\Lambda,V),\tilde\sigma))$. We obtain the
following result.

\begin{prop}\label{rangetr}
The range of the trace on $K_0(C^*(S(\Lambda,V),\tilde\sigma))$ is
\begin{equation}\label{trK0range}
[\tr](K_0(C^*(S(\Lambda,V),\tilde\sigma)))=\Z + \Z \theta(\theta'-\theta)^{-1}.
\end{equation}
\end{prop}

\proof
Since $X_\epsilon$ is a 3-manifold, when we expand the terms in the
cohomological formula \eqref{indrangetr} as
$$ \hat A(X_\epsilon)=1 - \frac{1}{24} p_1(X_\epsilon)+\cdots $$
$$ e^{\omega_\epsilon} = 1 + \omega_\epsilon + \frac{1}{2} \omega_\epsilon^2 +\cdots $$
$$ Ch(g)= -\frac{1}{6} \Tr(\beta(g)) + \frac{1}{5!} \Tr(\beta^3(g)) +
\cdots , $$
only the terms of the wedge product $\hat A(X_\epsilon) e^\omega
Ch(g)$ that give differential forms of order up to 3 can contribute 
nontrivially. 

Thus, we obtain the terms
\begin{equation}\label{intterms}
\frac{1}{(2\pi)^2} \int_{X_\epsilon} \left(\frac{-1}{6}
\Tr(\beta(g))\wedge \omega + \frac{1}{5!} \Tr(\beta(g)^3) \right). 
\end{equation}

The term 
$$ \frac{1}{(2\pi)^2} \int_{X_\epsilon} \frac{1}{5!} \Tr(\beta(g)^3) =
\frac{1}{(2\pi)^2} \int_{X_\epsilon} Ch(g) $$
is the term one would find in the untwisted case, and it gives the
untwisted odd Chern character.

For the remaining term
$$ \frac{1}{(2\pi)^2}\frac{-1}{6} \int_{X_\epsilon} 
\Tr(\beta(g))\wedge \omega_\epsilon $$
the range as $[g]$ varies in $K_1(C(X_\epsilon))$ is given by
$\Z R(\omega)$, where $R(\omega)$ is the range of the linear form 
$$ T_\omega: [g]\mapsto \frac{1}{(2\pi)^2}\frac{-1}{6} \int_{X_\epsilon} 
\Tr(\beta(g))\wedge \omega_\epsilon \in \R . $$
First notice that, with the notation
$Ch_1(g)=\frac{-1}{6}\Tr(\beta(g))$, we have
$$ \int_{C} Ch_1(g) =2\pi i \deg(g|_C) \in 2\pi i \Z, $$
for $C\in H_1(X_\epsilon,\Z)$. Thus, we obtain
$$ \frac{1}{2\pi i} \int_{X_\epsilon} Ch_1(g)\wedge PD(C) \in \Z, $$
for $PD(C)\in H^2(X_\epsilon,\Z)\hookrightarrow H^2(X_\epsilon,\R)$. 
Now consider the explicit description of the 2-form $\omega_\epsilon$
given in Lemma \ref{lemzetaomega} above. We can write $$\omega_\epsilon
= 2\pi i \theta(\theta'-\theta)^{-1} \bar\omega_\epsilon,$$ 
where $\bar\omega_\epsilon \in H^2(X_\epsilon,\Z)$ is given by
$$ 
\bar\omega_\epsilon (v,w) =A_\epsilon^k(\eta) \wedge \lambda,
$$
for $v=((0,0),(\lambda,k))$ and $w=((0,0),(\eta,r))$.
Thus, we see that we can write
$$ \frac{1}{(2\pi)^2} \int_{X_\epsilon} Ch_1(g)\wedge \omega_\epsilon
=-\frac{1}{2\pi i} \theta(\theta'-\theta)^{-1} \,\, 
\int_{X_\epsilon} Ch_1(g)\wedge \bar\omega_\epsilon $$
$$ = -\frac{1}{2\pi i}\theta(\theta'-\theta)^{-1} \,\,  
\int_{PD(\bar\omega_\epsilon)} Ch_1(g) =
\theta(\theta'-\theta)^{-1} \deg(g|_{PD(\bar\omega_\epsilon)}) 
\in \theta(\theta'-\theta)^{-1} \Z. $$
\endproof

\section{Isospectral deformations and spectral triples}\label{Sisospec}

In noncommutative geometry, the analog of Riemannian structures is
provided by the formalism of spectral triples \cite{Co-S3}. A spectral
triple on a noncommutative space $\cA$ (where $\cA$ is a
$C^*$-algebra) consists of the data $(\cA_\infty,\cH,D)$ of a dense
involutive subalgebra $\cA_\infty$, a representation $\pi:\cA\to \cB(\cH)$ as
bounded operators on a Hilbert space $\cH$ and a self-adjoint operator
$D$ on $\cH$, with compact resolvent, satisfying the compatibility
condition
\begin{equation}\label{boundcomm}
[D,\pi(a)] \in \cB(\cH), \ \ \ \ \forall a\in \cA_\infty.
\end{equation}

In particular, in the commutative case, to a Riemannian $spin$-manifold
$X$ one can associate a canonical spectral triple
$(C^\infty(X),L^2(X,S),\dirac)$. A reconstruction theorem
\cite{RenVar} shows that a spectral triple where the algebra is
abelian, which satisfies a list of axioms, is the canonical spectral
triple of a Riemannian $spin$-manifold.

In our case, we have a spectral triple associated to the 3-manifold
$X_\epsilon$, where the spinor bundle is a complex 2-plane bundle
and the Dirac operator can be written in the form
\begin{equation}\label{DiracX}
\dirac_{X_\epsilon} = c(dt) \frac{\partial}{\partial t} + c(e^t dx)
\frac{\partial}{\partial x} + c(e^{-t} dy) \frac{\partial}{\partial y}, 
\end{equation}
where $\{ dt, e^t dx, e^{-t} dy \}$ is the basis of the cotangent
bundle of $S(\R^2,\R,\epsilon)=\R^2\rtimes_\epsilon \R$ and
$c(\omega)$ denotes Clifford multiplication by the 1-form $\omega$.

More explicitly, \eqref{DiracX} is of the form
\begin{equation}\label{Diracepsilon}
\dirac_{X_\epsilon}= \frac{\partial}{\partial t} \sigma_0 +
e^t \frac{\partial}{\partial x} \sigma_1 + 
e^{-t} \frac{\partial}{\partial y} \sigma_2 =\left(\begin{matrix}
\frac{\partial}{\partial t} & e^{-t} \frac{\partial}{\partial y} -i
e^t \frac{\partial}{\partial x}  \\ e^{-t} \frac{\partial}{\partial y}
+i e^t \frac{\partial}{\partial x} & -\frac{\partial}{\partial t}
\end{matrix}\right), 
\end{equation}
where $\sigma_i$, for $i=0,1,2$, are the Pauli matrices. 

Our purpose here is to show that this commutative spectral triple can
be deformed isospectrally to a spectral triple for the noncommutative
tori $\bT_{\Lambda,i}$. 

\subsection{The Connes--Landi isospectral deformations}\label{Scola}

We consider the problem from the point of
view of the Connes--Landi isospectral deformations \cite{CoLa}. This
provides a general procedure to deform commutative spectral triples
to noncommutative ones isospectrally, for manifolds with isometric
torus actions.

We recall briefly the construction of isospectral deformations, in
a version that is best adapted to our setting.

Suppose given a spectral triple $(C^\infty(X),L^2(X,S),\dirac_X)$
associated to a compact Riemannian $spin$-manifold $X$. Assume that
the manifold $X$ has an action of a torus $T^2$ by isometries,
$T^2 \subset {\rm Isom}(X)$. Then one considers a noncommutative
algebra $\cA_\theta$, depending on a real parameter $\theta\in \R$,
which is obtained by decomposing the operators $\pi(f)\in \cB(\cH)$,
for $f\in C^\infty(X)$ and $\cH=L^2(X,S)$ according to their weighted
components 
\begin{equation}\label{fnm}
\pi(f)= \sum_{n,m\in \Z} \pi(f_{n,m}), 
\end{equation}
where
\begin{equation}\label{pifnm}
\alpha_{\tau}(\pi(f_{n,m})) = e^{2\pi i (n\tau_1 + m\tau_2)}\,
\pi(f_{n,m}), \ \ \ \forall \tau=(\tau_1,\tau_2)\in T^2,
\end{equation}
for
\begin{equation}\label{alphatau}
\alpha_\tau(T)=U(\tau)T U(\tau)^*, \ \ \ \forall T\in \cB(\cH),
\,\forall \tau\in T^2,
\end{equation}
with $U(\tau)$ the unitary transformations implementing the
$T^2$-action on $\cH=L^2(X,S)$ by
$$ U(\tau)\psi(x)=\psi(\tau^{-1}(x)). $$
Let $L_1$ and $L_2$ denote the infinitesimal generators of the action 
\begin{equation}\label{UtauL}
U(\tau)=\exp(2\pi i \tau L)=\exp(2\pi i(\tau_1 L_1 + \tau_2 L_2)). 
\end{equation}

We consider then the subalgebra of $\cB(\cH)$ generated by the
operators of the form 
\begin{equation}\label{deformxi12}
\pi_{\xi_1,\xi_2}(f) =\sum_{n,m} \pi(f_{n,m})  e^{-2\pi i (\xi_1 n
L_2+\xi_2 m L_1)},
\end{equation}
where $\xi_1$ and $\xi_2$ are two real parameters.

\begin{lem}\label{defNCtorus}
For homogeneous operators $\pi(f)_{n,m}$ define the deformed product
\begin{equation}\label{deformprod}
f_{n,m} *_{\xi_1,\xi_2} h_{k,r} :=e^{-2\pi i(\xi_1 nr+\xi_2 mk)} f_{n,m}
h_{k,r}.
\end{equation}
The product of operators of the form \eqref{deformxi12} satisfies
$$ \pi_{\xi_1,\xi_2}(f_{n,m})\pi_{\xi_1,\xi_2}(h_{k,r})
= \pi_{\xi_1,\xi_2}(f_{n,m} *_{\xi_1,\xi_2} h_{k,r}). $$
\end{lem}

\proof
One checks directly that the 
operator product $\pi(f_{\xi_1,\xi_2})\pi(h_{\xi_1,\xi_2})$ is 
given in components by
$$ \pi(f_{n,m}) *_{\xi_1,\xi_2} \pi(h_{k,r}) =
e^{-2\pi i(\xi_1 nr+\xi_2 mk)} \pi(f_{n,m}) \pi(h_{k,r}). $$
\endproof

One can recognize in \eqref{deformprod} the 
convolution product of the twisted group 
$C^*$-algebra $C^*(\Z^2,\sigma)$ with the cocycle
$$ \sigma((n,m),(k,r))=\exp(-2\pi i(\xi_1 nr+\xi_2 mk)). $$

As shown in \cite{CoLa}, the operators \eqref{deformxi12} have
bounded commutators with the Dirac operator. In fact, 
since $T^2$ acts by isometries, the Dirac operator satisfies
$$ U(\tau) D U(\tau)^* = D, $$
\ie it is of bidegree $(0,0)$. Thus, one sees that the commutators
$$ \begin{array}{rl}
[D,\pi_{\xi_1,\xi_2}(f)] = & \sum_{n,m} 
[D,\pi(f)_{n,m}e^{-2\pi i (\xi_1 n
L_2+\xi_2 m L_1)}] \\[2mm] = &
 \sum_{n,m} [D,\pi(f)]_{n,m}  e^{-2\pi i (\xi_1 n
L_2+\xi_2 m L_1)} ,
\end{array} $$
which is still a bounded operator on $\cH$. 

We consider in particular the case where $\xi_2 =u/2 =-\xi_1$.
We denote by $\cA_u =C^\infty(X)_u$ the deformed algebra, that is,
the algebra generated by the \eqref{deformxi12}.
The deformed spectral triple is given by the data
$(\cA_u,L^2(X,S),\dirac_X)$. 

\subsection{Noncommutative solvmanifolds}\label{Sncsolv}

We apply the procedure described above to obtain an isospectral
deformation of the solvmanifold $X_\epsilon$, which corresponds to
deforming the fiber tori to noncommutative tori.

The canonical spectral triple for $X_\epsilon$ consists of the data
$(C^\infty(X_\epsilon),L^2(X_\epsilon,S),\dirac_{X_\epsilon})$, with
the Dirac operator of the form \eqref{DiracX}. 

There is a torus action on $X_\epsilon$ by
isometries, which consists of translations along the fibers of the
fibration $T^2 \to X_\epsilon \to S^1$. This acts on spinors by
unitaries 
\begin{equation}\label{Utauxi}
U(\tau) \psi((x,y),t)=\psi((x+e^t\tau_1,y+e^{-t}\tau_2),t), 
\end{equation}
for $\tau\in T^2=\R^2/\Lambda$ and $(x,y)\in
T^2_t=\R^2/\Lambda_t$, the fiber over $t\in S^1$, with
$(e^t\lambda_1,e^{-t}\lambda_2)\in \Lambda_t$, for
$(\lambda_1,\lambda_2)\in \Lambda$.

The action clearly preserves the metric $dt^2+e^t dx^2 +
e^{-t} dy^2$ hence the Dirac operator \eqref{Diracepsilon}
satisfies
$$ U(\tau) \partial_{X_\epsilon} U(\tau)^* = \partial_{X_\epsilon}. $$
The infinitesimal generators of the action $\alpha_\tau$ are the
operators $2\pi L_1 =e^t \frac{\partial}{\partial x}$, $2\pi L_2
=e^{-t}\frac{\partial}{\partial y}$ with
$U(\tau)=\exp(2\pi i(\tau_1 L_1 + \tau_2 L_2))$. 

We introduce the following notation. We denote by $E_\lambda$,
for $\lambda\in \Lambda$, the function
\begin{equation}\label{Elambda}
 E_\lambda((x,y),t) :=
e^{2\pi i \langle\Theta_{-t}(x,y),\lambda\rangle} ,
\end{equation}
where, as above, $\Theta_{-t}(x,y)=(e^{-t}x, e^t y)$ and
$\langle (a,b),\lambda\rangle =a\lambda_1 + b\lambda_2$. We also
denote by $\Xi_u(\lambda, L_1,L_2)$ the operator
\begin{equation}\label{XiuL12}
\Xi_u(\lambda, L_1,L_2) :=
\exp\left(i\pi  \frac{u}{(\theta'-\theta)} \lambda \wedge
(L_1,L_2)\right),
\end{equation}
acting on $\cH=L^2(X_\epsilon,S)$.

\begin{prop}\label{Aualg}
The deformed algebra $C^\infty(X_\epsilon)_u$, for $u\in \R$, is the
$C^*$-subalgebra of $\cB(\cH)$, with $\cH=L^2(X_\epsilon,S)$ generated
by the operators of the form
\begin{equation}\label{udefoper}
\pi_u(f)= E_\lambda \Xi_u(\lambda, L_1,L_2) .
\end{equation}
\end{prop}

\proof The induced action $\alpha: T^2 \to \Aut(C^\infty(X_\epsilon))$
defined by $$\pi(\alpha_\tau(f))=U(\tau)\pi(f)U(\tau)^*$$ is of the form
$\alpha_\tau(f)((x,y),t)=f((x+e^t\tau_1,y+e^{-t}\tau_2),t)$. 

Thus, a homogeneous operator of bidegree 
$\lambda=(\lambda_1,\lambda_2)$ is in this
case a function $f_\lambda ((x,y),t)$ with the property that
\begin{equation}\label{fnmtau}
\alpha_\tau(f_\lambda) ((x,y),t) = e^{2\pi i 
(\lambda_1 \tau_1 + \lambda_2 \tau_2)} f_\lambda ((x,y),t).
\end{equation}
This condition is satisfied by functions of the form
\begin{equation}\label{exptnm}
f_\lambda ((x,y),t)=\exp(2\pi i\langle \Theta_{-t}(x,y),\lambda\rangle)=
\exp(2\pi i(e^{-t}\lambda_1 x+e^t \lambda_2 y)). 
\end{equation}

Under the change of variables 
\begin{equation}\label{changevar}
 \Z^2 \to \Lambda, \ \ \ (n,m)\mapsto \lambda=(n+m\theta,n+m\theta'), 
\end{equation}
the condition \eqref{fnmtau} corresponds to elements 
$f_{n,m}$ of bidegree $(n,m)$ for the corresponding 
action of $T^2=\R^2/\Z^2$.
Thus, using this change of coordinates to pass in \eqref{deformxi12}
from $\Z^2$ to $\Lambda$, 
we can see that elements of the deformed algebra of the form
\eqref{deformxi12} correspond to elements of the form
$$ \sum_\lambda a_\lambda \, E_\lambda\, \Xi_u(\lambda,L_1,L_2), $$
for $\xi_2 =u/2 =-\xi_1.$
\endproof

Set $2 u_\theta =u/(\theta'-\theta)$. The operators \eqref{udefoper}
act on spinors by
$$ \begin{array}{rl}
(\pi_u(f)\psi)((x,y),t)= & E_\lambda((x,y),t) 
(U((-\lambda_2 u_\theta ,\lambda_1 u_\theta))\psi)((x,y),t) \\[2mm]
= & e^{2\pi i\langle\Theta_{-t}(x,y),\lambda\rangle} 
\psi((x- e^t \lambda_2 u_\theta,y+ e^{-t}\lambda_1 u_\theta),t). 
\end{array} $$

\begin{prop}\label{NCtorusRep}
The operators
\begin{equation}\label{piRsigmalambda}
\pi(R^\sigma_\lambda) := E_\lambda \, 
 \Xi_u(\lambda,L_1,L_2)
\end{equation}
define a representation on $\cH=L^2(X_\epsilon,S)$ of the 
noncommutative torus $C^*(\Lambda,\sigma)$, with the cocycle
$$ \sigma(\lambda,\eta)=\exp(2\pi i u_\theta \, \lambda\wedge \eta). $$
\end{prop}

\proof
Notice that we have 
$$ U(-\lambda_2 u_\theta,\lambda_1 u_\theta) e^{2\pi i
\langle\Theta_{-t}(x,y),\eta\rangle} =
e^{2\pi i u_\theta \, \lambda\wedge \eta} e^{2\pi i
\langle\Theta_{-t}(x,y),\eta\rangle}. $$
Thus, we obtain
$$ e^{2\pi i \langle\Theta_{-t}(x,y),\lambda
\rangle} \Xi_u(\lambda,L_1,L_2) e^{2\pi i \langle\Theta_{-t}(x,y),\eta
\rangle} \Xi_u(\eta,L_1,L_2) = $$
$$ e^{2\pi i u_\theta \, \lambda\wedge \eta}
e^{2\pi i \langle\Theta_{-t}(x,y),\lambda\rangle}
e^{2\pi i\langle\Theta_{-t}(x,y),\eta\rangle}
\Xi_u(\lambda,L_1,L_2) \Xi_u(\eta,L_1,L_2)= $$
$$ e^{2\pi i u_\theta \, \lambda\wedge \eta} 
e^{2\pi i \langle\Theta_{-t}(x,y),\lambda+\eta\rangle}
\Xi_u(\lambda+\eta,L_1,L_2). $$
This shows that the operators $\pi(R^\sigma_\lambda)$ satisfy
the product rule 
$$ \pi(R^\sigma_\lambda)\pi(R^\sigma_\eta) =\sigma(\lambda,\eta) 
\pi(R^\sigma_{\lambda+\eta}), $$
for $\sigma(\lambda,\eta)=\exp(2\pi i u_\theta \, \lambda\wedge \eta)$,
which is the product rule of the twisted group algebra
$C^*(\Lambda,\sigma)$.
\endproof

We obtain in this way an isospectral noncommutative geometry given by
the finitely summable spectral triple
\begin{equation}\label{isosp3}
(C^\infty(X_\epsilon)_u,L^2(X_\epsilon,S),\dirac_{X_\epsilon}).
\end{equation}

\begin{cor}\label{isospecNCtor}
In the case $u=\theta$ and $u=\theta'$, the isospectral deformation 
\eqref{isosp3} defines a finitely summable spectral triple for 
the noncommutative tori $\bT_{\Lambda,i}$, with dense subalgebra
$\C(\Lambda,\sigma)$.
\end{cor}

\proof This is a direct consequence of Proposition \ref{NCtorusRep}
and the identifications of Corollary \ref{Lambdasigma} of the
$\bT_{\Lambda,i}$ with twisted group $C^*$-algebras $C^*(\Lambda,\sigma)$.
\endproof

The representation \eqref{piRsigmalambda} of $C^*(\Lambda,\sigma)$
extends to an action of $C^*(S(\Lambda,V),\tilde\sigma)$, as follows.
Let $U(k\log\epsilon)$ denote the unitary operator
\begin{equation}\label{Uk}
(U(k\log\epsilon)\psi)((x,y),t)=\psi(A_\epsilon^k(x,y),t)=
\psi((x,y),t-k\log\epsilon). 
\end{equation}

\begin{prop}\label{crossprodlaw}
The operators
\begin{equation}\label{piRsigmalambdak}
\pi(R^{\tilde\sigma}_{(\lambda,k)}) 
:= E_\lambda\, \Xi_u(\lambda,L_1,L_2)\, U(k\log\epsilon)
\end{equation}
define a representation on $\cH=L^2(X_\epsilon,S)$ of the twisted
group $C^*$-algebra $C^*(\Lambda\rtimes_\epsilon V,\tilde\sigma)$,
for the cocycle
$$ \tilde\sigma((\lambda,k),(\eta,r))=\exp\left(2\pi i
\frac{u_\theta}{2} \lambda \wedge A_\epsilon^k(\eta) \right). $$
\end{prop}

\proof We have the identities
$$ \begin{array}{l}
(U(k\log\epsilon) E_\eta) = E_{A_\epsilon^k(\eta)} , \\[2mm] 
\Xi_u(\lambda,L_1,L_2) E_{A_\epsilon^k(\eta)} =
e^{2\pi i u \lambda\wedge
A_\epsilon^k(\eta)} E_{A_\epsilon^k(\eta)} , \\[2mm]
 E_\lambda E_{A_\epsilon^k(\eta)} =E_{\lambda+A_\epsilon^k(\eta)} , \\[2mm]
 U(k\log\epsilon) \Xi_u(\eta,L_1,L_2) U(r\log\epsilon) =
\Xi_u(A_\epsilon^k(\eta),L_1,L_2) U((k+r)\log\epsilon). 
\end{array} $$ 
These combine to give the composition rule
$$ \begin{array}{c} E_\lambda \Xi_u(\lambda,L_1,L_2) U(k\log\epsilon) 
E_\eta \Xi_u(\eta,L_1,L_2) U(r\log\epsilon)= \\[2mm]
\tilde\sigma((\eta,r),(\lambda,k)) 
E_{\lambda+A_\epsilon^k(\eta)} \Xi_u(\lambda+
A_\epsilon^k(\eta),L_1,L_2) U((k+r)\log\epsilon). \end{array} $$
\endproof

\subsection{Unitary equivalences}\label{Sunitary}

We begin by reformulating the data described above in an
equivalent form by expanding in Fourier modes along the 
fiber tori as in \cite{ADS}.

Recall that the fiber over $t\in [0,\log\epsilon)$ is given 
by the torus $T^2_t=\R^2/\Lambda_t$, with
$\Lambda_t=\Theta_t(\Lambda)$. Thus, if we denote by $(x,y)$, as
above, the coordinates in $T^2_t$, we can write these as
$(x,y)=\Theta_t(a,b)$, with $(a,b)\in T^2=\R^2/\Lambda$,
the reference torus.

This means writing the spinors $\psi((x,y),t)$ in the form
\begin{equation}\label{spinFourier}
\sum_{\lambda} \psi_{\lambda} \,
e^{2\pi i\langle (a,b),\lambda \rangle} =\sum_{\lambda} \psi_{\lambda} \,
e^{2\pi i\langle \Theta_{-t}(x,y),\lambda \rangle} =\sum_{\lambda}
\psi_{\lambda} E_{\lambda}.
\end{equation}

The Dirac operator acts on $E_\lambda$ as
$$ \dirac_{X_\epsilon} E_\lambda =
(\frac{\partial}{\partial t}\sigma_0 + 2\pi i \lambda_1
\sigma_1 + 2\pi i \lambda_2) E_\lambda. $$
The operators $\pi(R^\sigma_\eta)$ act as 
\begin{equation}\label{actFourier}
E_\eta \Xi_u(\eta,L_1,L_2) E_\lambda =
e^{2\pi i u \eta\wedge\lambda} E_{\eta+\lambda}. 
\end{equation}
The commutators are bounded operators of the form
\begin{equation}\label{commFourier}
 [\dirac_{X_\epsilon},\pi(R^\sigma_\eta)]=
\left(\eta_1 \sigma_1 + \eta_2
\sigma_2\right) R^\sigma_\eta. 
\end{equation}
Thus, passing to Fourier modes in the fiber directions
gives a unitarily equivalent spectral triple for
the noncommutative tori $\bT_{\Lambda,i}$, with
\begin{equation}\label{unitary1}
\begin{array}{l}
\hat\dirac_{X_\epsilon} \psi_\lambda  = 
(\frac{\partial}{\partial t}\sigma_0 + 2\pi i \lambda_1
\sigma_1 + 2\pi i \lambda_2 \sigma_2) \psi_\lambda \\[2mm]
\hat\pi(R^\sigma_\eta) \psi_\lambda  = 
\sigma(\eta,\lambda)\, \psi_{\lambda+\eta}. 
\end{array}
\end{equation}

We then consider a second unitary equivalence, which, as in \cite{ADS}
adjusts for the possible signs of $\lambda_1$ and $\lambda_2$. Namely,
we define the following unitary operator on the Hilbert space of the
spinors $\psi_\lambda$. We set
\begin{equation}\label{U2}
\cU\, \psi_\lambda = \sigma_\lambda \, \psi_\lambda,
\end{equation}
where $\sigma_\lambda$ is a product of Pauli matrices, where
$\sigma_i$, for $i=1,2$ appears in the product if and only if
$\lambda_i <0$.
Then the Dirac operator transforms to the unitarily
equivalent operator
\begin{equation}\label{unitary2Dirac}
\cU \hat\dirac_{X_\epsilon} \cU^* = \sign(N(\lambda))
(\frac{\partial}{\partial t}\sigma_0 + 2\pi i |\lambda_1| 
\sigma_1 + 2\pi i |\lambda_2| \sigma_2) .
\end{equation}
The action of the $R^\sigma_\eta$ transform
correspondingly to the operators
\begin{equation}\label{unitary2R}
\cU \hat\pi(R^\sigma_\eta) \cU^* : \sigma_\lambda
\psi_\lambda \mapsto \sigma_{\lambda+\eta} 
\psi_{\lambda+\eta}. 
\end{equation}

We then perform the other unitary transformation used in \cite{ADS}.
To this purpose, let us fix a choice of a fundamental domain $\cF_V$
for the action of $V$ on the lattice $\Lambda$. By this choice of a fundamental
domain, we can write uniquely an element $\lambda\in \Lambda$ in the
form $\lambda=A_\epsilon^k(\mu)$, for a $\mu\in \cF_V$ and a $k\in\Z$.

For $\lambda=A_\epsilon^k(\mu)\neq 0$, 
consider then the time shift
\begin{equation}\label{Utimeshift}
\tilde\cU (\sigma_\lambda \psi_\lambda)(t) = \sigma_\lambda
\psi_\lambda (t-\log\frac{|\mu_1|}{|N(\mu)|^{1/2}}),
\end{equation}
so that we have
\begin{equation}\label{U3}
 \tilde\psi_\lambda := \tilde\cU (\sigma_\lambda \psi_\lambda)= 
\sigma_\lambda
\psi_{|N(\lambda)|^{1/2}(\sign(\lambda_1)\epsilon^k,
\sign(\lambda_2)\epsilon^{-k})}.
\end{equation} 

One obtains in this way a unitarily equivalent spectral triple for
$\bT_{\Lambda,i}$, with Dirac operator
\begin{equation}\label{D0non0}
\tilde \dirac =\tilde \dirac^{(0)} + \sum_{\mu \in (\Lambda 
\smallsetminus \{0 \})/V} \tilde \dirac^{(\mu)},
\end{equation} 
where 
\begin{equation}\label{U2dataD}
\begin{array}{l}
\tilde \dirac^{(\mu)}\,
\tilde\psi_{A_\epsilon^k(\mu)} = \\[2mm]
\sign(N(\mu)) |N(\mu)|^{1/2} \left( 
|N(\mu)|^{-1/2} \frac{\partial}{\partial t} \sigma_0 +2\pi i
\epsilon^k \sigma_1 + 
2\pi i \epsilon^{-k} \sigma_2 \right) \,  
\tilde\psi_{A_\epsilon^k(\mu)}
\end{array} 
\end{equation}
while the action of the $R^\sigma_\eta$ is by
$$  \tilde\pi(R^\sigma_\eta)\,
\tilde\psi_\lambda = \tilde\psi_{\lambda+\eta}. $$

As in \cite{ADS}, one can write 
the operator $\tilde \dirac^{(\mu)}$ as a product 
$$ \tilde \dirac^{(\mu)} =D_\mu B_\mu, $$ with
\begin{equation}\label{DB}
\begin{array}{l}
D_\mu\,  \tilde\psi_{A_\epsilon^k(\mu)} = \sign(N(\mu))
|N(\mu)|^{1/2}\, \tilde\psi_{A_\epsilon^k(\mu)} \\[2mm]
B_\mu\,  \tilde\psi_{A_\epsilon^k(\mu)} = \left( 
|N(\mu)|^{-1/2} \frac{\partial}{\partial t} \sigma_0 + 2\pi
i\epsilon^k \sigma_1 + 2\pi i \epsilon^{-k}\sigma_2 \right)
\tilde\psi_{A_\epsilon^k(\mu)} . 
\end{array}
\end{equation}

In the following, we relate the Dirac operator 
$\dirac_{X_\epsilon}$, its unitarily equivalent operators 
discussed here above, and the decomposition \eqref{DB}
to known differential operators on noncommutative tori. 

\subsection{Differential operators on noncommutative tori}\label{DiffNCtori}

Notice that the action on $\R^2$ of the 1-paramater subgroup of $\SL_2(\R)$ 
$$ \Theta_t =\left(\begin{matrix} e^t & 0 \\ 0 & e^{-t}
\end{matrix}\right) $$ 
has fixed point $(0,0)$, with stable manifold the axis $(0,y)$ and unstable
manifold the axis $(x,0)$. On the standard torus $\R^2/\Z^2$ with
coordinates $(s_1,s_2)$ with $(x,y)=(s_1+s_2\theta,s_1+s_2\theta')$
these two directions define the two Kronecker foliations
$s_1+s_2\theta$ and $s_1+s_2\theta'$ with conjugate slopes $\theta$
and $\theta'$. The points of the lattice $\Lambda$ determine on these
two foliations the points of the pseudolattices $\Z+\Z\theta$ and
$\Z+\Z\theta'$, which define the equivalence relation on the space of
leaves of the two Kronecker foliations, defining as quotients the
noncommutative tori $\bT_{\Lambda,i}$, $i=1,2$. 
The action of $\Theta_t$ is expanding along the
line $L_\theta =\{ s_1+s_2\theta \}$ and contracting along
$L_{\theta'} =\{ s_1+s_2\theta' \}$ and flows the other points of
$\R^2$ along hyperbola with asymptotes $L_{\theta'}$ and $L_\theta$. 

\smallskip

Thus, the operators $e^t \frac{\partial}{\partial x}$ and $e^{-t}
\frac{\partial}{\partial y}$ correspond to derivations along the leaf
direction of these two transverse Kronecker foliations. The factors
$e^t$ and $e^{-t}$ are the normalization factors that account for the
rescaling of the transverse measure due to the action of the flow
$\Theta_t$. In fact, consider for instance a small transversal of
length $\ell$ for the Kronecker foliation $L_\theta$, given by the
interval $T_\ell=\{ (x,y): x=1,\,  -\ell/2 < y < 
\ell/2 \}$. The flow $\Theta_t$ maps it to the transversal
$\Theta_t(T_\ell)=\{ (x=e^t, y): \, -e^{-t}\ell/2< y < e^{-t}\ell/2
\}$ of length $e^{-t}\ell$. Thus, the differentiation
$\frac{\partial}{\partial x}$ in the leaf direction of $L_\theta$
is weighted by the factor $e^t$ that normalizes the length of the
transversal and corrects for the scaling of the transverse measure. 

\smallskip

Consider then the terms $2\pi i \lambda_1 \sigma_1$ and
$2\pi i \lambda_2\sigma_2$ in the operator 
$$ \hat\dirac_{X_\epsilon}: 
\psi_\lambda \mapsto (\frac{\partial}{\partial t}\sigma_0 + 2\pi i
\lambda_1 \sigma_1 + 2\pi i \lambda_2\sigma_2)\,\psi_\lambda, $$ 
that we obtained after passing to Fourier modes on the
fiber tori $T^2_t$. These terms correspond, respectively, to 
the leafwise derivations $e^t \frac{\partial}{\partial x}$ and $e^{-t}
\frac{\partial}{\partial y}$. These can be expressed equivalently in
terms of the operators
$$  \delta_\theta : \psi_{n,m}\mapsto (n+m\theta)\, \psi_{n,m}, \ \ \
\text{ and } \ \ \  \delta_{\theta'}: \psi_{n,m}\mapsto (n+m\theta')\,
\psi_{n,m}, $$
so that the sum $ \lambda_1 \sigma_1+ \lambda_2\sigma_2$
acts as the operator
\begin{equation}\label{Dtheta}
 \Dirac_{\theta,\theta'}= 
\left(\begin{matrix} 0 & \delta_{\theta'}-i \delta_\theta \\
\delta_{\theta'} + i \delta_\theta & 0 \end{matrix}\right). 
\end{equation}
This gives the Dirac operator of a spectral triple on the
noncommutative tori $\bT_{\Lambda,i}$ with
$$ 
R^\sigma_{r,k} \psi_{n,m}= \sigma((r,k),(n,m)) \psi_{(n,m)+(r,k)}
$$
and 
$$ [ \Dirac_{\theta,\theta'},R^\sigma_{r,k}] = 
\left(\begin{matrix} 0 & (r+k\theta')-i (r+k\theta) \\
(r+k\theta') + i (r+k\theta) & 0 \end{matrix}\right)R^\sigma_{r,k}. $$
In the particular case where $\theta'=-\theta$, this 
agrees with the spectral triple for the first order signature 
operator on the noncommutative torus considered, for instance, 
in \cite{MRR}. The construction of
\cite{MRR} can be interpretaed as obtained by using the two transverse
Kronecker foliations $L_\theta$ and $L_{-\theta}$ and the associated 
leafwise derivations $\partial/\partial x$ and $\partial/\partial y$. 

We can consider here the same kind of unitary transformations that
we described earlier for $\hat\dirac_{X_\epsilon}$, applied to
the operator $\Dirac_{\theta,\theta'}$ of \eqref{Dtheta}. Let us denote
by $\Dirac_{\theta,\theta',0}$ the restriction of $\Dirac_{\theta,\theta'}$ to 
the complement of the zero modes $\psi_0$ (\ie $\lambda=0$).
We have, as in \eqref{D0non0},
\begin{equation}\label{Dtheta0non0}
\Dirac_{\theta,\theta',0} =\sum_{\mu \in (\Lambda \smallsetminus \{ 0 \})/V} 
\Dirac^\mu_{\theta,\theta'},
\end{equation}
with
$$ \Dirac^\mu_{\theta,\theta'} \, \psi_{A_\epsilon^k(\mu)} = 
(\lambda_1 \sigma_1+ \lambda_2\sigma_2) \, 
\psi_{A_\epsilon^k(\mu)}. $$
After the unitary
transformation $\tilde\cU\cU$ with $\cU$ as in \eqref{U2} and
$\tilde\cU$ as in \eqref{U3}, we obtain a unitarily equivalent 
operator
\begin{equation}\label{tildeDtheta}
\tilde \Dirac^\mu_{\theta,\theta'} \, \tilde \psi_{A_\epsilon^k(\mu)} =
\sign(N(\mu))\, |N(\mu)|^{1/2} \, (\epsilon^k \sigma_1 
+ \epsilon^{-k} \sigma_2) \, \tilde \psi_{A_\epsilon^k(\mu)}.
\end{equation} 
As before, we factor this as a product of the operators
\begin{equation}\label{DBprodtheta}
\tilde \Dirac^\mu_{\theta,\theta'}= D^\mu_\theta B_\theta,
\end{equation}
with
\begin{equation}\label{DBtheta}
\begin{array}{l}
D^\mu_\theta \, \tilde \psi_{A_\epsilon^k(\mu)} = \sign(N(\mu))\,
|N(\mu)|^{1/2} \, \tilde \psi_{A_\epsilon^k(\mu)} \\[2mm]
B_\theta \, \tilde \psi_{A_\epsilon^k(\mu)} = (\epsilon^k \sigma_1 
+ \epsilon^{-k} \sigma_2) \, \tilde \psi_{A_\epsilon^k(\mu)}.
\end{array}
\end{equation}

\section{Shimizu $L$-function and Lorentzian geometry}\label{Slorentz}

In this section we describe another way of relating the Shimizu $L$-function to
the geometry of noncommutative tori with real multiplication, by regarding the norms $N(\lambda)$,
for $\lambda\in \Lambda$, as defining the momenta of a {\em Lorentzian} rather
than Euclidean Dirac operator.

\smallskip

Instead of working with positive inner
product spaces, as in the case of Euclidean spectral triples, the Galois involution of the real 
quadratic field defines a natural choice of a ``Krein involution'' and the norm correspondingly
defines an indefinite quadratic form. One formulates in this way a notion of spectral triple
over a real quadratic field and with Lorentzian signature, using the relation between 
indefinite inner product spaces and the associated real Hilbert spaces. The main point that 
requires care is the fact that the Lorentzian Dirac operator has a noncompact group of symmetries, 
in our case given by the units of the real quadratic field, hence it fails to have compact resolvent
due to the presence of infinite multiplicities in the eigenvalues. We show that the multiplicities
can be resolved by transforming the triple via a Krein isometry that is an unbounded self-adjoint
operator in the associated real Hilbert space and defines a finitely summable associated Dirac
operator in the Euclidean signature. 

\medskip 

As above, we let $\Lambda$ be the lattice in $\R^2$ associated to a lattice $L\subset \K$ 
in a real quadratic field $\K=\Q(\sqrt{d})$ by the embeddings $\iota_i:\K \hookrightarrow \R$,
\begin{equation}\label{LambdaL}
\Lambda=\{ \lambda \in \R^2\, |\, \lambda=(\lambda_1,\lambda_2)=(\iota_1(\ell),\iota_2(\ell)), 
\, \ell\in L\}. 
\end{equation}
We denote, as above, by $V$ the group $V=\epsilon^\Z$ of units preserving $\Lambda$. 
We denote the action as above with $\lambda \mapsto A_\epsilon^k(\lambda)
=(\epsilon^k\lambda_1,\epsilon^{-k}\lambda_2)$.

For $x\in \K$, we denote by $x'=c(x)$ the image under the Galois involution of $\K$. 
For $\lambda=(\lambda_1,\lambda_2)\in \Lambda$, we have $\lambda_2=c(\lambda_1)$. 
The norm is given by $N(\lambda)=\lambda_1\lambda_2$, and $N(\epsilon)=\epsilon\epsilon'=1$.

We consider the quadratic form $N(\lambda)=\lambda_1 \lambda_2=
(n+m\theta)(n+m\theta')$ to be the analog of the wave operator $\Box=p_0^2 -p_1^2$.

Its Dirac factorization into linear first order operator is obtained by considering 
a linear operator of the form
\begin{equation}\label{Dlambda}
\cD_\lambda =\left(\begin{matrix} 0 & \cD_\lambda^+ \\ \cD_\lambda^- & 0
\end{matrix}\right) 
:=\left(\begin{matrix} 0 & \lambda_1 \\ \lambda_2 & 0
\end{matrix}\right), 
\end{equation}
whose square is $\cD^2_\lambda =\Box_\lambda$, with
\begin{equation}\label{boxlambda}
 \Box_\lambda =\left(\begin{matrix} N(\lambda) & 0 \\ 0 & N(\lambda)
\end{matrix}\right). 
\end{equation}

We assemble these modes to define an operator $\cD$ acting on
$\cH=\ell^2(\Lambda)\oplus\ell^2(\Lambda)$ by $\cD e_{\lambda,\pm} 
=\cD_\lambda e_{\lambda,\pm}$. This satisfies $\cD \gamma =-\gamma
\cD$ with respect to the $\Z/2\Z$-grading 
$$ \gamma =\left(\begin{matrix} 1 & 0 \\ 0 & -1 \end{matrix}\right). $$

Consider the algebra $C^*(\Lambda,\sigma)$ of the noncommutative torus
acting diagonally on $\cH$. The operator $\cD$ has bounded commutators
with the elements of the dense subalgebra $\C(\Lambda,\sigma)$ since
we have
\begin{equation}\label{boundcommsigma}
 [\cD,R^\sigma_\eta] e_{\lambda,\pm}= \sigma(\lambda,\eta) \,
\eta_{\mp} \,\, e_{\eta+\lambda,\pm}, 
\end{equation}
where we used the notation $\eta_+=\eta_1$ and $\eta_-=\eta_2$.

However, the other properties of $\cD$ differ significantly from what
one usually postulates for Dirac operators of spectral triples. 

First of all,  notice that $\cD$ is not self-adjoint. In fact, it is invariant with
respect to a different involution, defined for operators with
coefficients in the real quadratic field $\K$, namely
$$ \cD = c(\cD^t), $$
where $\cD^t=(\cD^t_\lambda)$ denotes the transpose and $c(\cD)$ denotes the
effect of the Galois involution $c: x\mapsto x'$ of $\K$ applied 
to the coefficients of $\cD$. In this arithmetic context, it is
natural to require this property instead of self-adjointness. 

A more serious problem, however, comes from the fact that the operator $\cD$ has
infinite multiplicities, hence it is very far from having the compact
resolvent property of spectral triple. This is a typical problem
one encounters in trying to extend the formalism of spectral triples
from the Euclidean to the Lorentzian context, because of the presence
of non-compact symmetry groups for Lorentzian manifold. Here the
non-compact symmetry group is given by the units in $V =\epsilon^\Z$. 

\subsection{Arithmetic Krein spaces}\label{Skrein}

It is well known that, when one replaces Euclidean geometry by
Lorentzian geometry, the notion of Hilbert space is replaced by the
notion of a Krein space (\cf \eg \cite{Bognar}). The version we consider here
is slightly different from the usual one, since we want to be able to work over
the real quadratic field $\K$ instead of passing directly to complex numbers.

\begin{defn}\label{Kreinprod}
Let $c:\K\to \K$ denote the Galois involution $c:x \mapsto x'$ of the
real quadratic field. Let $\cV$ be a $\K$-vector space. We say that a map $T:\cV\to \cV$ 
is $c$-linear if it satisfies $T(av+bw)=c(a) T(v)+ c(b) T(w)$.

A Lorentzian pairing on a $\K$-vector space $\cV$ is a 
non-degenerate $\K$-valued pairing
$$ (\cdot,\cdot): \cV\times \cV \to \K $$
which is $c$-linear in the first variable and linear in the second. 
\end{defn}

We can then introduce the analog of the notion of a Krein space, 
in this arithmetic context. 

\begin{defn}\label{KreinK}
A Krein space over a real quadratic field $\K$ (or $\K$-Krein space) is a 
$\K$-vector space $\cV$ endowed with a Lorentzian pairing 
$(\cdot,\cdot)$ as in 
Definition \ref{Kreinprod}, and a $c$-linear involution 
$\kappa:\cV\to \cV$, such that the pairing $(\kappa\cdot,\cdot)$ has the following properties:
\begin{enumerate} 
\item $(\kappa\cdot,\cdot)=c(\cdot,\kappa\cdot)$
\item For all $v\neq 0$ in $\cV$, the elements
$(\kappa v,v)\in \K$ are totally positive.
\end{enumerate}
\end{defn}

When properties (1) and (2) of Definition \ref{KreinK} holds, we say that
$(\kappa\cdot,\cdot)$ defines a positive definite inner product. We have 
a corresponding notion of Krein adjoint as follows.

\begin{defn}\label{Kreinadjdef}
Given a $\K$-linear operator $T$ on a $\K$-Krein space $\cV$, 
the Krein adjiont $T^\dag$ is the adjoint in the Lorentzian pairing 
$(\cdot,\cdot)$, 
\begin{equation}\label{kreinadjdef}
(v,T w) = (T^\dag v, w).
\end{equation}
\end{defn}

The $c$-linear involution $\kappa$ of Definition \ref{KreinK} corresponds to a
Wick rotation from Lorentzian to Euclidean signature. 
The Krein adjoint satisfies $T^\dag=\kappa T^* \kappa$, where
$T^*$ is the adjoint in the inner product $\langle \cdot,\cdot\rangle=
(\kappa\cdot,\cdot)$. 

Given a $\K$-Krein space $\cV$, there are two naturally associated real Hilbert spaces, 
obtained by considering the real vector spaces 
\begin{equation}\label{VRi}
\cV_{\R,i} := \cV \otimes_{\iota_i(\K)} \R,
\end{equation}
obtained by tensoring $\cV$ with $\R$ using either one of the two embeddings
$\iota_i:\K \hookrightarrow \R$ of the real quadratic field.

\begin{lem}\label{VRinnprod}
The pairing
\begin{equation}\label{innprod12}
\langle v,w \rangle = \frac{1}{2}  \iota_1\left((\kappa v,w) +(v,\kappa w)\right) =
\frac{1}{2}  \iota_2\left((\kappa v,w) +(v,\kappa w)\right)
\end{equation}
induced on $\cV_{\R,i}$ by the Lorentzian pairing $(\cdot,\cdot)$ on $\cV$ defines
a real valued positive definite inner product.
\end{lem}

\proof We know that $(v,\kappa w)=c(\kappa v, w)$. Thus, we have 
$$
\langle v, w \rangle =\frac{1}{2} \left( \iota_1(\kappa v,w) + \iota_2(\kappa v,w) \right).
$$
We can extend this pairing by $\R$-linearity to define a bilinear form on $\cV_{\R,i}$.
By the assumption that for $v\neq 0$ the $(\kappa v,v)$ are totally positive elements of $\K$ 
we obtain that \eqref{innprod12} defines a positive-definite inner product. 
\endproof

in the following, we still denote by $\cV_{\R,i}$ the Hilbert space completion 
obtained in this way.

\subsection{Lorentzian spectral triples over real quadratic fields}\label{SfieldS3}

It is not uncommon to make use of Krein spaces to extend the formalism of 
spectral triples to Lorentzian geometry, see \eg \cite{Stro}. 
Here we follow a similar viewpoint, adapted to the arithmetic setting of real quadratic fields.

For a $\K$-linear operator $T$ acting on  a $\K$-Krein space $\cV$, we define
$\bM_i(T)\geq -\infty$, for $i=1,2$, as
\begin{equation}\label{Mopnorm}
\bM_i(T):= \inf_{(v,v)=1} \iota_i (T v, T v).
\end{equation}

We introduce the following preliminary notion of a $\K$-triple, which we then 
refine by additional properties providing the analog of a spectral triple.

\begin{defn}\label{preS3}
A Krein $\K$-triple consists of data $(\cA,\cV,\cD)$ with the following properties. 
\begin{enumerate}
\item $\cA$ is an involutive algebra over the real quadratic field $\K$. 
\item $\cV$ is a $\K$-Krein space with non-degenerate $\K$-bilinear form $(\cdot,\cdot)$.
\item The algebra $\cA$ acts on $\cV$ via a 
representation $\pi: \cA \to \End_\K(\cV)$,  
with the involution of $\cA$ realized by the Krein adjoint $\pi(a^*)=\pi(a)^\dag$.
\item The operators $\pi(a)$, for $a\in\cA$, satisfy 
\begin{equation}\label{bMa}
\bM_i(a)>-\infty.
\end{equation}
\item The operator $\cD$ is a densely defined $\K$-linear operator on $\cV$, 
which is Krein-self-adjoint , $\cD^\dag=\cD$.  
\item The commutators $C_a:=[\cD,a]$ satisfy
\begin{equation}\label{Kreinboundcomm1}
 \bM_i(C_a) > -\infty,  \ \ \ \forall a\in \cA.
\end{equation} 
\end{enumerate}
\end{defn}

We then define Lorentzian $\K$-spectral triples in the following way.

\begin{defn}\label{KreinS3}
A Krein $\K$-triple $(\cA,\cV,\cD)$ as in Definition \ref{preS3} is a finitely summable 
Lorentzian $\K$-spectral triple if the following holds. 
\begin{enumerate}
\item There exists a densely defined $\K$-linear operator $U:\cV \to \cV$ with 
$(Uv,Uv)=(v,v)$, for all $v\in {\rm Dom}(U)$ and $U^\dag =U^{-1}$, with the property that
\begin{equation}\label{D2symmU}
U^\dag \cD U = \cD. 
\end{equation}
\item The commutators $C_{a,U}:=[\cD_U,\pi_U(a)]$, with $\pi_U(a)=U^\dag \pi(a) U$, 
satisfy the condition
\begin{equation}\label{Kreinboundcomm}
\bM_i(C_{a,U}) > -\infty, \ \ \ \forall a\in \cA.
\end{equation} 
\item The operator $U$ is an unbounded self-adjoint operator, $U=U^*$ on the associated real 
Hilbert space $\cV_{\R,i}$ with the inner product $\langle\cdot,\cdot\rangle$ of \eqref{innprod12}.  
\item The triple $(\cA,\cV,\cD)$ is $p$-summable for 
$p\in \R^*_+$ if
\begin{equation}\label{finsumKrein}
\sum_n \left| \langle  U e_n, |\cD^2| U e_n \rangle \right|^{-s/2}   <\infty, 
\ \ \ \ \forall s \geq p,
\end{equation}
where $e_n$ is an orthonormal basis for the complement of the zero modes of 
the operator $|\cD^2|$ in the real Hilbert space $\cV_{\R,i}$. 
\end{enumerate}
\end{defn}

Notice that in Krein spaces isometries are not necessarily bounded 
operators (see \cite{Bognar} \S VI), so the $U$ is only densely defined in general. 
The definition given here is different from the notions of Lorentzian spectral 
triples currently developed in the literature. The differences stem mainly 
from our need to work over a finite extension of $\Q$ instead of $\C$ and 
to resolve the infinite multiplicity of the eigenvalues. We also require the 
weaker property \eqref{bMa}, \eqref{Kreinboundcomm1} and \eqref{Kreinboundcomm}, 
instead of requiring continuity in the operator norm in the associted Hilbert space. 
These conditions will become more transparent in our main example below. 

\subsection{Arithmetic twisted group algebras}\label{Sartwist}

We consider the $\K$-vector space $\cV_\Lambda$ spanned by
the basis elements $e_\lambda$ with $\lambda\in \Lambda$, endowed with the
pairing
\begin{equation}\label{bilform}
(v,w):= \sum_\lambda c(a_\lambda) b_\lambda,
\end{equation}
for $v=\sum_\lambda a_\lambda e_\lambda$ and  $w=\sum_\lambda b_\lambda e_\lambda$, 
and with $c: x\mapsto x'$ the Galois involution of $\K$.

\begin{lem}\label{KreinHilbLambda}
The space $\cV_\Lambda$ with the pairing \eqref{bilform} is a $\K$-Krein space.
\end{lem}

\proof Clearly the pairing \eqref{bilform} is a Lorentzian pairing in the sense of
Definition \ref{Kreinprod}.
Let $\kappa:\cV_\Lambda\to \cV_\Lambda$ be given by the Galois involution
$$ \kappa (v) = \sum_\lambda c(a_\lambda) e_\lambda. $$
Then the pairing $\langle v,w\rangle =(\kappa v,w)=c(v,\kappa w)$ is a 
positive definite inner product, as in Definition \ref{KreinK}. In fact, we have
$$ \iota_1 \langle v,v\rangle =\sum_\lambda \iota_1(a_\lambda^2) \geq 0, \ \ \  
\iota_2 \langle v,v\rangle =\sum_\lambda \iota_2(a_\lambda^2)\geq 0. $$
\endproof

We consider on $\cV_\Lambda$ the action of the group ring $\K[\Lambda]$,
given by $R_\lambda e_\eta = e_{\lambda +\eta}$. 

\begin{lem}\label{boundedopsR}
The operators $R_\lambda$ acting on $\cV_\Lambda$ satisfy
$$ \bM_i(R_\lambda) > -\infty. $$
Moreover, the operators $R_\lambda$ define bounded operators 
in the associated real Hilbert spaces 
$\cV_{\Lambda,\R,i}=\cV_\Lambda\otimes_{\iota_i(\K)}\R$. 
\end{lem}

\proof The operators $R_\lambda$ are Krein isometries, 
and $(R_\lambda v,R_\lambda v)=(v,v)$ implies that $\bM_i(R_\lambda)=1$. 
The operators $R_\lambda$ act by $e_\eta \mapsto e_{\lambda+\eta}$ on the 
associated Hilbert spaces, hence they define bounded (unitary) operators.
\endproof

Now we want to introduce, in this setting of $\K$-Krein spaces, an analog of 
the twisted group ring $\C(\Lambda,\sigma)$ (the noncommutative torus) 
we have been working with in the complex case. 

\begin{lem}\label{Kcocycle}
Suppose given $\omega \in \K^*$ with $N(\omega)=\omega \omega'=1$. Then the expression
\begin{equation}\label{sigmaomega}
\varpi(\lambda,\eta)= \omega^{(n,m)\wedge (r,k)},
\end{equation}
for $\lambda=(n+m\theta,n+m\theta')$ and $\eta=(r+k\theta,r+k\theta')$, defines a
$\K^*$-valued group 2-cocycle $\sigma$ on $\Lambda$.
\end{lem}

\proof The argument is the same as in the complex case. It suffices to show that 
the cocycle condition holds. \endproof

\begin{defn}\label{Ktwistgr}
The twisted group ring $\K(\Lambda,\sigma)$ is the unital involutive $\K$-algebra 
generated by elements $R^\varpi_\lambda$ with
the product
\begin{equation}\label{KRsigma}
R^\varpi_\lambda R^\varpi_\eta = \varpi(\lambda,\eta)
R^\varpi_{\lambda+\eta} = \omega^{(n,m)\wedge (r,k)}
R^\varpi_{\lambda+\eta},
\end{equation}
for $\lambda=(n+m\theta,n+m\theta')$ and $\eta=(r+k\theta,r+k\theta')$, and the involution
$(R^\varpi_\lambda)^*=R^\varpi_{-\lambda}$.
\end{defn}

The twisted group ring $\K(\Lambda,\sigma)$ also acts on $\cV_\Lambda$ by
\begin{equation}\label{varpiact}
R^\varpi_\lambda e_\eta = \varpi(\eta,\lambda) e_{\lambda +\eta}.
\end{equation}

\begin{lem}\label{varpiactlem}
The operators $R^\varpi_\lambda$, acting as in \eqref{varpiact}, satisfy 
$\bM_i(R^\varpi_\lambda)>-\infty$. 
\end{lem}

\proof The action \eqref{varpiact} preserves the
Lorentzian pairing $(\cdot,\cdot)$ on $\cV_\Lambda$ since 
$$ (R^\varpi_\lambda e_\eta, R^\varpi_\lambda e_\zeta)=
c(\omega(\eta,\lambda))\omega(\zeta,\lambda)\delta_{\eta,\zeta} = 
N(\omega(\eta,\lambda))\delta_{\eta,\zeta}=(e_\eta,e_\zeta), $$
since $N(\omega(\eta,\lambda))=1$. The condition
$(R^\varpi v,R^\varpi v)=(v,v)$ implies $\bM_i(R^\varpi_\lambda)=1$.
\endproof

However, notice that, while the action of $\K[\Lambda]$ 
extends to an action by bounded operators on
the associated real Hilbert spaces $\cV_\K\otimes_{\iota_i(\K)}\R$,
the induced action on $\cV_\R$ of the twisted group ring 
$\K(\Lambda,\varpi)$ is by the unbounded operators 
\begin{equation}\label{RsigmaVR}
R^\varpi_\lambda\, e_{\eta,\pm} = A_\omega^{(r,k) \wedge (n,m)}\,
e_{\lambda+\eta,\pm}, 
\end{equation}
with
\begin{equation}\label{Aomega}
A_\omega = \left(\begin{matrix} \iota_1(\omega) & 0 \\ 0 &
\iota_2(\omega)  \end{matrix}\right) \in \SL_2(\R).
\end{equation}

As in the complex case, we can also consider the group ring 
$\K[S(\Lambda,V)]$ for $S(\Lambda,V)=\Lambda\rtimes V$. 
The cocycle \eqref{sigmaomega} extends to a cocycle on the cross product by setting 
\begin{equation}\label{tildesigmaomega}
\tilde\varpi((\lambda,k),(\eta,r))= \varpi(\lambda,A_\epsilon^k(\eta))=
\omega^{(n,m)\wedge (u,v)\varphi_\epsilon^k},
\end{equation}
for $\lambda=(n+m\theta,n+m\theta')$ and $\eta=(u+v\theta,u+v\theta')$, with $n,m,u,v\in\Z$.

\begin{defn}\label{twistomegaS}
The twisted group ring $\K(S(\Lambda,V),\tilde\varpi)$ is the 
unitary involutive $\K$-algebra with generators $R^{\tilde\varpi}_{\lambda,k}$ satisfying
$$ R^{\tilde\varpi}_{\lambda,k}R^{\tilde\varpi}_{\eta,r} =\tilde\varpi((\lambda,k),(\eta,r))
R^{\tilde\varpi}_{\lambda+A_\epsilon^k(\eta),k+r}, $$
with the involution $(R^{\tilde\varpi}_{\lambda,k})^*=
R^{\tilde\varpi}_{-A_\epsilon^{-k}(\lambda),-k}$.
\end{defn}

\subsection{Lorentzian Dirac operator}\label{SlorentzD}

On the $\K$-Krein space $\cV_\Lambda \oplus \cV_\Lambda$ 
we consider the densely defined $\K$-linear operator 
\begin{equation}\label{Dlambdac}
\cD_\K e_{\lambda,\pm} = \cD_{\K,\lambda}e_{\lambda,\pm}  =
\left(\begin{matrix} 0 & \cD_\lambda^+ \\ \cD_\lambda^- & 0
\end{matrix}\right) e_{\lambda,\pm} 
:=\left(\begin{matrix} 0 & \ell \\ c(\ell)  & 0
\end{matrix}\right) e_{\lambda,\pm} , 
\end{equation}
where we write $\lambda\in\Lambda$ as $\lambda = (\iota_1(\ell),\iota_2(\ell))$ 
with $\ell\in L\subset \K$, as in \eqref{LambdaL}.

The operator $\cD_\K$ of \eqref{Dlambdac} induces on the real Hilbert space 
$\cV_{\R,i}\oplus \cV_{\R,i}$ the $\R$-linear operators
\begin{equation}\label{DRi}
\cD_\lambda =\left(\begin{matrix} 0 & \lambda_1 \\ \lambda_2  & 0
\end{matrix}\right)  \ \ \ \text{ and } \ \ \ c(\cD_\lambda)=
\left(\begin{matrix} 0 & \lambda_2 \\ \lambda_1  & 0
\end{matrix}\right),
\end{equation}
respectively. 
This recovers the Lorentzian Dirac operator described in \eqref{Dlambda} above.

\begin{lem}\label{AVD1}
The data $(\K(\Lambda,\varpi),\cV_\Lambda\oplus\cV_\Lambda,\cD_\K)$ 
define a Krein $\K$-triple in the sense of Definition \ref{preS3}.
\end{lem}

\proof Properties (1)--(4) of Definition \ref{preS3} 
follow from Lemma \ref{KreinHilbLambda}, Lemma \ref{varpiactlem}, 
and the fact that the Krein adjoint $(R^\varpi_\lambda)^\dag =
R^\varpi_{-\lambda}=(R^\varpi_\lambda)^{-1}$. 
Property (5) follows directly from \eqref{Dlambdac}, since
$$ \cD_\K^\dag =c(\cD_{\K,\lambda}^t)=
\left(\begin{matrix} 0 & c(\lambda_2) \\ c(\lambda_1)  & 0
\end{matrix}\right) =
\left(\begin{matrix} 0 & \lambda_1 \\ \lambda_2  & 0 \end{matrix}\right) =\cD_\K. $$ 
We then need to prove (6), namely that the commutators $[\cD_\K, R^\varpi_\lambda]$ satisfy
$$ \bM_i([\cD_\K, R^\varpi_\lambda]) > -\infty. $$
We have
$$ [\cD_\K, R^\varpi_\lambda] e_{\eta,\pm} = \varpi(\eta,\lambda) 
\left(\begin{matrix} 0 & (\lambda_1+\eta_1)-\eta_1 \\  (\lambda_2+\eta_2)-\eta_2  & 0
\end{matrix}\right) e_{\lambda+ \eta,\pm}. $$
Thus, we have
$$ ( [\cD_\K, R^\varpi_\lambda] v, [\cD_\K, R^\varpi_\lambda] v) =N(\lambda) (v,v), $$
from which the result follows.
\endproof

Suppose given a choice of a fundamental domain $\cF_V$ for the action of 
$V=\epsilon^\Z$ on $\Lambda$. Let $\rho(\lambda)\in\Z$ denote the unique integer such that
$\lambda=A_\epsilon^{\rho(\lambda)}(\mu)$, with $\mu\in \cF_V$.

Consider the $\K$-linear operator on $\cV_\Lambda\oplus \cV_\Lambda$ defined by
\begin{equation}\label{Tupsilon} 
T_\epsilon e_{\lambda,\pm} := \left(\begin{matrix} \epsilon^{\rho(\lambda)} & 0 \\ 0 &
 \epsilon^{-\rho(\lambda)} \end{matrix}\right) e_{\lambda,\pm} .
\end{equation}
Consider also the involution $J:\cV\to \cV$ defined by setting
$$ J e_{\lambda,\pm} = e_{J(\lambda),\pm}, $$
where $J(\lambda) =A_\epsilon^{-k}(\mu)$ for 
$\lambda=A_\epsilon^k(\mu)$ with $\mu\in \cF_V$ and $k\in\Z$. 
This satisfies $J^2=1$ and $J^\dag =J$.

We set $U_\epsilon = T_\epsilon J$, with
\begin{equation}\label{Uepsilon}
U_\epsilon e_{\lambda,\pm} = \left(\begin{matrix} \epsilon^{-\rho(\lambda)} & 0 \\ 0 &
 \epsilon^{\rho(\lambda)} \end{matrix}\right) e_{J(\lambda),\pm}.
\end{equation}

We now show that the data of Lemma \ref{AVD1} satisfy the properties a Lorentzian $\K$-spectral triple.

\begin{prop}\label{AVD2}
The data $(\K(\Lambda,\varpi),\cV_\Lambda\oplus\cV_\Lambda,\cD_\K)$ define a 
Lorentzian $\K$-spectral triple,  as in Definition \ref{KreinS3}.
\end{prop}

\proof  The $T_\epsilon$ are Krein isometries, since 
$$ (T_\epsilon e_{\lambda,\pm}, T_\epsilon e_{\lambda,\pm})=  N(\epsilon^{\rho(\lambda)})  
(e_{\lambda,\pm},e_{\lambda,\pm}) = (e_{\lambda,\pm},e_{\lambda,\pm}). $$
They satisfy $T^\dag_\epsilon = T_{\epsilon^{-1}} = T_\epsilon^{-1}$. Thus we have
$U_\epsilon^\dag = J^\dag T_\epsilon^\dag = J T_\epsilon^{-1} = U_\epsilon^{-1}$. 
This is also a Krein isometry since both $T_\epsilon$ and $J$ are, with ${\rm Dom}(U_\epsilon)=
{\rm Dom}(T_\epsilon)$, since $J$ is bounded.

The operator $U_\epsilon$ is a symmetry of the Dirac operator, namely we have 
\begin{equation}\label{DKupsilon1}
\cD_{\K,\epsilon}:= U_\epsilon^\dag \cD_\K U_\epsilon =\cD_\K. 
\end{equation}

In fact, we have
\begin{equation}\label{DKomega}
T_\epsilon^\dag \cD_\K T_\epsilon e_{\lambda,\pm} = 
\left(\begin{matrix} 0& \epsilon^{-2\rho(\lambda)} \lambda_1 \\ 
\epsilon^{2\rho(\lambda)} \lambda_2 & 0 \end{matrix}\right) e_{\lambda,\pm}.
\end{equation}
Since we have $\lambda=A_\epsilon^{\rho(\lambda)}(\mu)$ with $\mu\in \cF_V$, we can write 
the above equivalently as 
$$ T_\epsilon^\dag  \cD_\K T_\epsilon e_{\lambda,\pm} = 
\left(\begin{matrix} 0& \epsilon^{-\rho(\lambda)} \mu_1 \\ 
\epsilon^{\rho(\lambda)} \mu_2 & 0 \end{matrix}\right) e_{\lambda,\pm}. $$
Thus, we have
$$ J^\dag T_\epsilon^\dag  \cD_\K T_\epsilon J e_{\lambda,\pm} =
\left(\begin{matrix} 0& \epsilon^{\rho(\lambda)} \mu_1 \\ 
\epsilon^{-\rho(\lambda)} \mu_2 & 0 \end{matrix}\right) e_{\lambda,\pm} $$
$$ 
= \left(\begin{matrix} 0& \lambda_1 \\ 
\lambda_2 & 0 \end{matrix}\right) e_{\lambda,\pm} =\cD_\K e_{\lambda,\pm}.
$$
This proves property (1) of Definition \ref{KreinS3}.

(2) follows from Lemma \ref{AVD1} and the fact that $U_\epsilon$ is a Krein isometry, since
$$ (U^\dag_\epsilon [\cD_\K,R^\varpi_\lambda] U_\epsilon  v, U^\dag_\epsilon 
 [\cD_\K,R^\varpi_\lambda] U_\epsilon  v)=
([\cD_\K,R^\varpi_\lambda] U_\epsilon  v, [\cD_\K,R^\varpi_\lambda] U_\epsilon  v) $$ 
$$ = N(\lambda) (U_\epsilon v,U_\epsilon v)=
N(\lambda)(v,v). $$

(3) The adjoint $U^*_\epsilon$ in the associated Hilbert space inner product 
$\langle \cdot,\cdot\rangle$ of \eqref{innprod12}
satisfies $U^*_\epsilon =U_\epsilon$ since
$$ U^*_\epsilon e_{\lambda,\pm}= (\kappa U^\dag_\epsilon \kappa) e_{\lambda,\pm}= c(U^\dag_\epsilon) e_{\lambda,\pm} $$ 
$$ = \left(\begin{matrix} c(\epsilon^{\rho(\lambda)}) & 0 \\ 0 &
 c(\epsilon^{-\rho(\lambda)}) \end{matrix}\right) e_{J(\lambda),\pm} =  
\left(\begin{matrix} \epsilon^{-\rho(\lambda)} & 0 \\ 0 &
 \epsilon^{\rho(\lambda)} \end{matrix}\right) e_{J(\lambda),\pm} = U_\epsilon e_{\lambda,\pm} . $$
 
Consider then the operator $|\cD_\K^2|$ acting on the associated real Hilbert space by
$$ |\cD_\K^2| e_{\lambda,\pm} =
\left(\begin{matrix} |N(\lambda)| & 0 \\ 
0 & |N(\lambda)| \end{matrix} \right) e_{\lambda,\pm}. $$
We restrict $|\cD_\K^2|$ to the orthogonal complement of the 
zero modes, \ie on the span of the $e_{\lambda,\pm}$ with 
$\lambda\neq 0$. We then obtain
\begin{equation}\label{Zetasum}
\sum_{\lambda\neq 0}  \left| \langle U_\epsilon e_{\lambda,\pm},  |\cD_\K^2| 
U_\epsilon e_{\lambda,\pm} \rangle \right|^{-s/2}  = 
\sum_{\lambda\neq 0} (\epsilon^{2\rho(\lambda)}
+\epsilon^{-2\rho(\lambda)})^{-s/2} |N(\lambda)|^{-s/2}.
\end{equation}
This can be written equivalently as  
\begin{equation}\label{Zetasum2}
  \sum_{k\in\Z} (\epsilon^{2k}+\epsilon^{-2k})^{-s/2} 
 \sum_{\mu \in (\Lambda\smallsetminus \{ 0 \})/V}
 |N(\mu)|^{-s/2}, 
\end{equation}
using the unique decomposition $\lambda=A_\epsilon^k(\mu)$, 
for $k\in\Z$ and $\mu\in \cF_V$, associated to the choice of the fundamental domain. 
Thus, we see that the finite summability condition holds.
\endproof

\smallskip

\begin{defn}\label{etaLdef}
The eta function of a Lorentzian $\K$-spectral triple is the function
\begin{equation}\label{etaLorentz}
\eta_\cD(s):= \sum_n  \sign(\langle U e_n, \cD^2 U e_n \rangle)
\left|\langle U e_n, |\cD^2| U e_n \rangle \right|^{-s/2} ,
\end{equation}
where the sum is over an orthonormal basis for the complement 
of the zero modes of $|\cD^2|$ in the Hilbert space $\cV_{\R,i}$.
\end{defn}

The following result relates the Shimizu $L$-function to the Lorentz $\K$-spectral triple.

\begin{cor}\label{ShimLorentz}
The eta function for the Lorentz $\K$-spectral triple of Proposition \ref{AVD2} is of the form
\begin{equation}\label{etaDK}
\eta_{\cD_\K}(s) =  L(\Lambda,V,s/2) Z_\epsilon(s/2),
\end{equation}
where $L(\Lambda,V,s)$ is the Shimizu $L$-function 
and $Z_\epsilon(s/2) =\sum_{k\in\Z} (\epsilon^{2k}+\epsilon^{-2k})^{-s/2}$.
\end{cor}

\proof The argument is the same as in Proposition \ref{AVD2}. We have
$$  \sum_{\lambda\neq 0} \sign(\langle U_\epsilon e_{\lambda,\pm}, \cD_\K^2 U_\epsilon e_{\lambda\pm} \rangle)
\left| \langle U_\epsilon e_{\lambda,\pm},  |\cD_\K^2| U_\epsilon e_{\lambda,\pm} \rangle\right|^{-s/2} $$
$$ = \sum_{\lambda\neq 0} \sign(N(\lambda)) (\epsilon^{2\rho(\lambda)}
+\epsilon^{-2\rho(\lambda)})^{-s/2} |N(\lambda)|^{-s/2} $$
 $$ = \sum_{k\in\Z} (\epsilon^{2k}+\epsilon^{-2k})^{-s/2} \sum_{\mu \in 
 (\Lambda\smallsetminus \{ 0 \})/V}
 \sign(N(\mu)) |N(\mu)|^{-s/2}. $$
The result then follows since $L(\Lambda,V,s)=\sum_{\mu \in (\Lambda\smallsetminus \{ 0 \})/V}
 \sign(N(\mu)) |N(\mu)|^{-s}$.
\endproof

\subsection{Eta function and 3-dimensional geometry}\label{Seta3d}

The zeta and eta functions we obtained in Proposition \ref{AVD2} 
and Corollary \ref{ShimLorentz} for the Lorentzian spectral geometry 
are closely related to those one can obtain from the spectral geometry 
of the 3-dimensional solvmanifold $X_\epsilon$ and the signature 
operator on the noncommutative torus. 

We have seen in \S \ref{DiffNCtori} above that the Dirac 
operator $\dirac_{X_\epsilon}$ on the 3-dimensional 
solvmanifold $X_\epsilon$ can be related to the signature 
operator $\Dirac_{\theta,\theta'}$ of \eqref{Dtheta} on the 
noncommutative torus $\bT_{\Lambda,i}$. 
Up to a unitary equivalence, we have written in \eqref{DBprodtheta} 
the operator $\Dirac_{\theta,\theta'}$ in terms of the operators 
$\tilde \Dirac^\mu_{\theta,\theta'}= D^\mu_\theta B_\theta$, 
with $\mu \in (\Lambda \smallsetminus \{ 0 \})/V$, defined as in  \eqref{DBtheta}. 

The zeta an eta functions for the operator 
$\tilde \Dirac^\mu_{\theta,\theta'}$ have the following form.

\begin{lem}\label{zetaetaDtheta}
The operator $\tilde \Dirac_{\theta,\theta'}$ has an 
associated zeta function of the form
\begin{equation}\label{zetaDtheta}
\zeta_{\tilde \Dirac_{\theta,\theta'}}(s) = 
2 Z_\epsilon(s/2) \sum_{\mu \in (\Lambda \smallsetminus \{ 0 \})/V} |N(\mu)|^{-s/2}.
\end{equation}
The eta function of $\tilde \Dirac_{\theta,\theta'}$ vanishes 
due to the symmetry in the spectrum. However, the restriction 
$\tilde \Dirac_{\theta,\theta'}^+$ of $\tilde \Dirac_{\theta,\theta'}$ 
to the subspace $\cH^+$ of the positive modes of the operator $B_\theta$ 
has a nonvanishing eta function of the form
\begin{equation}\label{etaDthetaplus}
\eta_{\tilde \Dirac_{\theta,\theta'}^+}(s)    
= L(\Lambda,V,s/2) Z_\epsilon(s/2) = \eta_{\cD_\K}(s).
\end{equation}
\end{lem}

\proof The operator
$$ B_\theta =\left( \begin{matrix} 0 & \epsilon^{-k} -i \epsilon^k \\
\epsilon^{-k} + i \epsilon^k & 0  \end{matrix} \right) $$
has spectrum 
$$ \Spec(B_\theta) =\{ \pm (\epsilon^{2k}+\epsilon^{-2k})^{1/2} \} $$
which is symmetric around zero. Thus, for the zeta function we have
$$ \zeta_{\tilde \Dirac_{\theta,\theta'}}(s) = \sum_{\mu \in 
(\Lambda \smallsetminus \{ 0 \})/V} |N(\mu)|^{-s/2} \, 
2 \sum_{k\in\Z} (\epsilon^{2k}+\epsilon^{-2k})^{-s/2}, $$
while the eta function vanishes. 

One can restrict the spectral triple for
the noncommutative torus $\bT_{\Lambda,i}$ 
to the subspace $\cH^+$ of the positive modes
of the operator $B_\theta$, since the action of the $R^\sigma_\eta$
preserve this decomposition. The new Dirac operator
$\tilde \Dirac^+_{\theta,\theta'}$  is then given by the restriction of
$\tilde \Dirac_{\theta,\theta'}$ to $\cH^+$. It has a corresponding 
decomposition
$$ \tilde \Dirac^+_{\theta,\theta',0} 
=\sum_{\mu \in (\Lambda \smallsetminus \{0
\})/V} D^\mu_\theta B_\theta^+, $$
where $\Spec(B_\theta^+)=\{ (\epsilon^{2k}+\epsilon^{-2k})^{1/2} \}$.
Thus, in this case one obtains
$$ \eta_{\tilde \Dirac^+_{\theta,\theta'}}(s) = 
L(\Lambda,V,s/2)\, \zeta_{B_\theta^+}(s), $$
where
$$ \zeta_{B_\theta^+}(s) =Z_\epsilon(s/2)=\sum_k
(\epsilon^{2k}+\epsilon^{-2k})^{-s/2}. $$
\endproof

\subsection{The residue}\label{Sresidue}

The special value $L(\Lambda,V,0)$ of the Shimizu $L$-function can be extracted 
from the eta function $\eta_{\tilde \Dirac^+_{\theta,\theta'}}(s)$ 
in the following way.

\begin{cor}\label{residue}
The function $\eta_{\tilde \Dirac^+_{\theta,\theta'}}(s)$ has a pole
of order one at $s=0$ with
\begin{equation}\label{Res}
\Res_{s=0}\,\, \eta_{\tilde \Dirac^+_{\theta,\theta'}}(s) =
\frac{L(\Lambda,V,0)}{\log \epsilon}.
\end{equation}
\end{cor}

\proof Consider the function  
$$ Z_\epsilon(s):= \sum_{k\in\Z} (\epsilon^{2k} + \epsilon^{-2k})^{-s}. $$
It suffices to show that it has a simple pole at $s=0$ with residue 
\begin{equation}\label{ZezRes}
\Res_{s=0} Z_\epsilon(s) = \frac{1}{\log \epsilon}.
\end{equation}
One writes
$$ \Gamma(s) Z_\epsilon (s) = \int_0^\infty g_\epsilon(t) \, t^{s-1} dt,  $$
where
\begin{equation}\label{fbt}
 g_\epsilon(t)=  \left(\sum_{k\in\Z} e^{-(\epsilon^{2k}+\epsilon^{-2k})t} \right) 
\end{equation}
for $t>0$. 
The function $g_\epsilon(t)$ satisfies
$$ g_\epsilon(t)= -e^{-2t} + 2 h_\epsilon(t) - 2 \sum_{k=0}^\infty e^{-\epsilon^{2k} t}
(1- e^{-\epsilon^{-2k}t}), $$
where 
$$ h_\epsilon(t)=\sum_{k=0}^\infty e^{-\epsilon^{2k} t}. $$
We can estimate $-e^{-2t} = -1 + O(t)$ when $t\to 0$ and
$(1-e^{-\epsilon^{-2k}t}) = O(\epsilon^{-2k} t)$, uniformly.
Notice that 
$$ h_\epsilon(t)-h_\epsilon (\epsilon^2 t)=e^{-t}=\sum_{r=0}^\infty \frac{(-1)^r}{r!} t^r $$ 
hence 
\begin{equation}\label{htasympt}
 h_\epsilon(t)=\frac{1}{2\log \epsilon} \log(1/t) + C - \sum_{r=0}^\infty 
\frac{(-1)^r}{r! (\epsilon^{2r}-1)} t^r . 
\end{equation}
Thus, the function $\Gamma(s)Z_\epsilon(s)$ has a double
pole at $s=0$ and simple poles at $s\in \Z_{<0}$.
Thus, the function $Z_\epsilon(s)$ has a simple pole at zero with 
residue $1/\log \epsilon$.
\endproof


\begin{thebibliography}{99}

\bibitem{At} M.F.~Atiyah, {\em Elliptic operators, discrete groups and
von Neumann algebras}.  
Colloque ``Analyse et Topologie" en l'Honneur de Henri Cartan (Orsay,
1974),  pp. 43--72. Asterisque, No. 32-33, Soc. Math. France, Paris, 1976. 

\bibitem{ADS} M.F.~Atiyah, H.~Donnelly, I.M.~Singer, {\em 
Eta invariants, signature defects of cusps, and values of $L$-functions},
Annals of Mathematics 118 (1983) 131--177.

\bibitem{APS} M.F.~Atiyah, V.K.~Patodi, I.M.~Singer, {\em Spectral 
asymmetry and Riemannian geometry. III}.  Math. Proc. Cambridge 
Philos. Soc.  79  (1976), no. 1, 71--99. 

\bibitem{BC2} P.~Baum, A.~Connes, {\em Chern character for discrete groups}.  
A fete of topology,  163--232, Academic Press, Boston, MA, 1988. 

\bibitem{Bel1} J.~Bellissard, {\em Noncommutative geometry and 
quantum Hall effect}.  Proceedings of the International Congress 
of Mathematicians, Vol. 1, 2 (Z\"urich, 1994),  1238--1246, 
Birkh\"auser, 1995. 

\bibitem{Bel3} J.~Bellissard, {\em Gap labelling theorems for
Schr\"odinger operators},  From number theory to physics (Les Houches,
1989),  538--630, Springer, Berlin, 1992.

\bibitem{Bognar} J.~Bogn\'ar, {\em Indefinite inner product spaces},
Springer, 1974.

\bibitem{ChaEch} J.~Chabert, S.~Echterhoff, {\em Permanence properties
of the Baum--Connes conjecture}, Doc. Math. 6 (2001) 127--183.

\bibitem{CoCR} A.~Connes, {\em $C^*$-alg\`ebres et g\'eom\'etrie 
diff\'erentielle}.  
C. R. Acad. Sci. Paris S\'er. A-B 290 (1980), no. 13, A599--A604.

\bibitem{Co-thom} A.~Connes, {\em An analogue of the Thom isomorphism 
for crossed products of a $C^*$-algebra by an action of $\R$}.
Adv. in Math. 39 (1981), no. 1, 31--55.

\bibitem{Co-S3} A.~Connes, {\em Geometry from the spectral point of
view}.  Lett. Math. Phys.  34  (1995),  no. 3, 203--238. 

\bibitem{CoLa} A.~Connes, G.~Landi, {\em Noncommutative manifolds, 
the instanton algebra and isospectral deformations}. 
Comm. Math. Phys. 221 (2001), no. 1, 141--159. 

\bibitem{ELPW} S.~Echterhoff, W.~L\"uck, N.C.~Phillips, S.~Walters,
{\em The structure of crossed products of irrational rotation algebras
by finite subgroups of $SL_2(\Z)$}, preprint SFB 478, Heft 435, 2006.

\bibitem{Getz} E.~Getzler, {\em The odd Chern character in cyclic homology 
and spectral flow}. Topology 32 (1993), no. 3, 489--507. 

\bibitem{Hirz} F.~Hirzebruch, {\em Hilbert modular surfaces}.  
Enseignement Math. (2)  19  (1973), 183--281.

\bibitem{hofs} D.G.~Hofstadter, {\em Energy levels and wave functions
of Bloch electrons in rational or irrational magnetic
field}. Phys. Rev. B14 (1976) 2239--2249.

\bibitem{Man} Yu.I.~Manin, {\em 
Real multiplication and noncommutative geometry (ein Alterstraum)}.  
The legacy of Niels Henrik Abel,  685--727, Springer, Berlin, 2004.

\bibitem{MarMat1} M.~Marcolli, V.~Mathai, {\em Twisted index theory on good
orbifolds, I: noncommutative Bloch theory}, Commun. Contemp. Math. Vol1 (1999)
N.4, 553--587.

\bibitem{MarMat2} M.~Marcolli, V.~Mathai, {\em Twisted index theory on good
orbifolds, II: fractional quantum numbers}, Commun. Math. Phys. 217 (2001)
55--87.

\bibitem{MarMat3} M.~Marcolli, V.~Mathai, {\em Towards the fractional
quantum Hall effect, a noncommutative geometry perspective}, 
Noncommutative Geometry and Number Theory, pp.235--262. Vieweg Verlag,
2006.

\bibitem{MRR} R.~Matthes, O.~Richter, G.~Rudolph, {\em Spectral triples and 
differential calculi related to the Kronecker foliation}.  
J. Geom. Phys.  46  (2003),  no. 1, 48--73.

\bibitem{MiFo} A.~Mishchenko, A.~Fomenko, {\em The index of elliptic 
operators over $C^*$-algebras}. Math. USSR Izv. 15 (1980) 87--112.

\bibitem{Muller} W.~M\"uller, {\em Signature defects of cusps of
Hilbert modular varieties and values of $L$-series at $s=1$}.  
J. Differential Geom.  20  (1984),  no. 1, 55--119. 

\bibitem{PaRae} J.A.~Packer, I.~Raeburn, {\em On the structure of
twisted group $C^*$-algebras}, Trans. AMS 334 (1992) 685--718.

\bibitem{RenVar} A.~Rennie, J.~Varilly, {\em Reconstruction of 
manifolds in noncommutative geometry}, arXiv:math/0610418.

\bibitem{Rie} M.A.~Rieffel, {\em $C^*$-algebras associated with 
irrational rotations}.  Pacific J. Math.  93  (1981), no. 2,
415--429. 

\bibitem{Ros} J.~Rosenberg, {\em $C^*$-algebras, positive scalar
curvature, and the Novikov conjecture III}, Topology 25 (1986)
319--336. 

\bibitem{Stro} A.~Strohmaier, {\em On noncommutative and
pseudo-Riemannian geometry}.  J. Geom. Phys.  56  (2006),  
no. 2, 175--195.

\bibitem{Sunada} T.~Sunada, {\em A discrete analogue of periodic 
magnetic Schr\"odinger operators}. Geometry of the spectrum (Seattle, 
WA, 1993),  283--299, Contemp. Math., 173, Amer. Math. Soc., 
Providence, RI, 1994. 


\end{thebibliography}
\end{document}